\documentclass[3p]{elsarticle}
\journal{Journal of Computational Physics}

\biboptions{sort&compress} 

\usepackage[bookmarks=true,colorlinks=true,linkcolor=blue]{hyperref}
\usepackage{amsmath,amsfonts,amssymb,mathabx,bm}
\usepackage{graphicx,esint,ulem}
\usepackage{calc}
\usepackage{geometry}
\usepackage[usenames,dvipsnames,svgnames,table]{xcolor}
\usepackage{stmaryrd}   

\usepackage{color}

\usepackage{algorithm,algpseudocode}
\usepackage{algorithmicx}
\algrenewcommand\alglinenumber[1]{\footnotesize #1:} 

\usepackage{tikz}

\usepackage{tcolorbox}
\usepackage[english]{babel}
\usepackage{fullpage,setspace,float}
\usepackage{tabu}
\usepackage{multirow}
\usepackage{mdframed}

\newcommand{\scst}{\beta} 
\newcommand{\Mtyp}{M_{\rm{typ}}}

\newcommand{\Nd}[1]{\mathcal{N}_{\delta}(#1)}
\newcommand{\Ndt}[1]{\mathcal{N}_{\Delta t}(#1)}
\newcommand{\floor}[1]{\left\lfloor #1 \right\rfloor}

\newcommand{\Oh}{\mathcal O}
\newcommand{\sigvec}{\bm{\sigma}}
\newcommand{\gvec}{\bm{g}}
\newcommand{\ones}{\bm{1}\bm{1}^{T}}
\newcommand{\err}[2]{\varepsilon_{#1}^{#2}}
\newcommand{\con}[2]{\Delta_{#1}^{#2}}
\newcommand{\tcon}[2]{\tilde\Delta_{#1}^{#2}}
\newcommand{\errvec}[1]{\bm\varepsilon^{#1}}

\newcommand{\pt}{\partial_t}
\newcommand{\px}{\partial_x}

\newcommand{\dt}{\Delta t}

\newcommand{\beq}{\begin{equation}}
\newcommand{\eeq}{\end{equation}}
\newcommand{\beqs}{\begin{equation*}}
\newcommand{\eeqs}{\end{equation*}}

\newcommand{\Nt}{N_t}

\newcommand{\abs}[1]{\left\vert #1\right\vert}

\newcommand{\norm}[1]{\Vert #1\Vert}

\newcommand{\uin}{u_{\text{inc}}}
\newcommand{\utot}{u_{\text{tot}}}
\newcommand{\jump}[1]{\left[ #1 \right]}

\DeclareMathOperator{\dist}{dist}
\DeclareMathOperator{\supp}{supp}
\DeclareMathOperator{\sinc}{sinc}

\newcommand{\R}{\mathbb{R}}
\newcommand{\C}{\mathbb{C}}
\newcommand{\Z}{{\mathbb Z}}
\newcommand{\al}{\alpha}
\newcommand{\tal}{\tilde\alpha}
\newcommand{\si}{\sigma}
\newcommand{\eps}{\epsilon}
\newcommand{\om}{\omega}
\newcommand{\bigO}{\mathcal{O}}
\newcommand{\qqquad}{\qquad\qquad}

\newcommand{\hatbump}{\widehat{\phi'}}          
\newcommand{\LLR}{{L^2(\R)}}            
\newcommand{\LLT}{{L^2(0,T)}}
\newcommand{\bi}{\begin{itemize}}
\newcommand{\ei}{\end{itemize}}
\newcommand{\mmin}{{m_{\text{min}}}}
\newcommand{\mmax}{{m_{\text{max}}}}
\newcommand{\ntyp}{{M_{\text{typ}}}}        

\newtheorem{theorem}{Theorem}

\newtheorem{definition}{Definition}
\newenvironment{proof}[1][Proof]{\begin{trivlist}
\item[\hskip \labelsep {\bfseries #1.}]}{\end{trivlist}}
\newtheorem{thm}{Theorem}[section]

\newtheorem{lem}{Lemma}[section]

\newtheorem{remark}{Remark}[section]

\begin{document}

\def\ccm{Center for Computational Mathematics, Flatiron Institute, Simons Foundation,
  New York, New York 10010}

\def\nyu{Courant Institute of Mathematical Sciences,
  New York University, New York, New York 10012}

\begin{frontmatter}
 \title{A fast algorithm for the wave equation using time-windowed Fourier projection}

 \author{Nour G. Al Hassanieh\fnref{ccm,nyu}}
\address[ccm]{\ccm}
\address[nyu]{\nyu}
 \ead{nalhassanieh@flatironinstitute.org}

\author{Alex H. Barnett\fnref{ccm}}
\ead{abarnett@flatironinstitute.org}

\author{Leslie Greengard\fnref{ccm,nyu}}
\ead{lgreengard@flatironinstitute.org}

\begin{abstract}
  We introduce a new
arbitrarily high-order
  method for the rapid evaluation of hyperbolic
potentials (space-time integrals involving the Green's 
function for the scalar wave equation).
With $M$ points in the spatial discretization and $N_t$ time steps of size $\dt$,
a naive implementation would require
$\bigO(M^2N_t^2)$ work in dimensions where the weak Huygens' principle applies.
We avoid this all-to-all interaction
using a smoothly windowed decomposition into
a local part, treated directly, plus a history part, approximated by a
$N_F$-term Fourier series.
In one dimension, our method requires 
$\bigO\left((M + N_F \log N_F)N_t\right)$ work,
with ${N_F} = \bigO(1/\dt)$, by exploiting the
non-uniform fast Fourier transform.
We demonstrate the method's performance for time-domain scattering problems involving
a large number $M$ of springs (point scatterers) attached
to a vibrating string at arbitrary locations, with either periodic or free-space boundary conditions.
We typically achieve 10-digit accuracy, and include tests for $M$ up to a million.
\end{abstract}


\end{frontmatter}

\section{Introduction} 

In this paper, we consider the scalar wave equation 
in one space dimension,
posed either in free space ($x \in \mathbb{R}$) or in a periodic box
($x \in B$, with $B = [-\pi,\pi]$).
This equation governs the propagation of waves on a
string, where the function of interest is 
the displacement of the string in the vertical direction.
It also governs linear acoustics and (through a slight generalization)
time-dependent electromagnetics.
The numerical solution of such wave problems, especially
the two- and three-dimensional (3D) cases, and especially at high frequency
(many wavelengths across the domain),
are extremely useful in engineering applications. We are particularly interested
in advancing the integral equation approach
\cite{fmmbook}.
The acceleration techniques we present do generalize to higher dimensions, but
the application to scattering demands other complications such as spatial quadratures.
Thus we chose to present the ideas in the setting of a high-order accurate solution
of a 1D model problem,
but one which retains the challenges of high frequencies and multiple scattering.

\begin{figure}[th]  \label{fig:spring}   
\centering
\includegraphics[width=.9\textwidth]{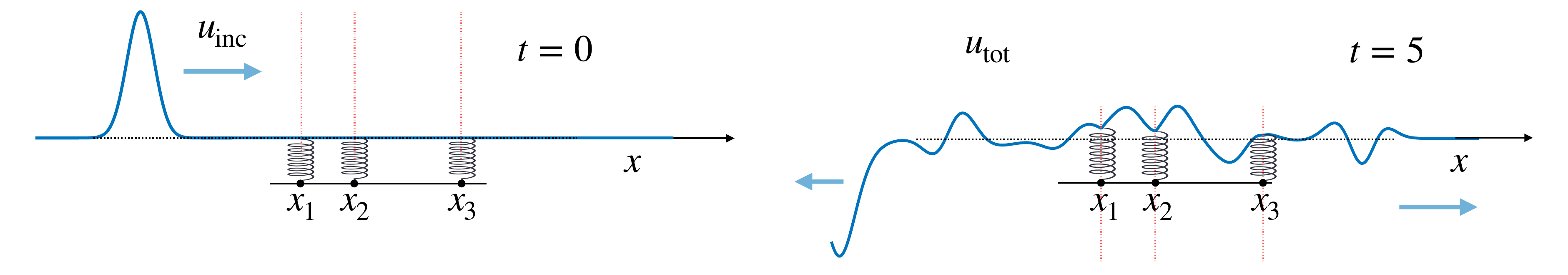}
\caption{(Left) A wideband pulse propagates rightwards along an infinite
string, and impinges on a region where the string is attached to a set
of springs, i.e., point-like scatterers.
(The equilibrium position aligns with the quiescent string.)
(Right) The total displacement is plotted at a later time after the initial
pulse has traversed the springs. Note that the displacement is continuous
but not smooth. (The jump in slope is 
most visible at the locations of the first two springs.) }
\end{figure}

While the pure initial value problem for the wave equation admits an
analytic solution \cite[Ch.~4 \& 10]{Guenther2012}, most problems of interest involve boundary and/or
interface conditions. We consider point-like scatterers
for 1D waves. At each point $x_j$ in a collection
$\Gamma = \{ x_j\}_{j = 1}^{M} \subset B$,
the string is attached vertically to a fixed reference by a spring of
nondimensionalized spring constant $\beta_j>0$,
as sketched in Fig.~\ref{fig:spring}.
(The natural length of the spring is chosen such that its force
is zero when the string is quiescent.)
A given rightward-traveling solution of the homogeneous wave equation 
\beq
\uin(x,t) = f(x-t)
\label{uinc}
\eeq
then impinges on this structure.
We will assume that $f$ is smooth and that $\Gamma \cap \supp f = \emptyset$.
In the standard scattering theory setup,
the total (physical) displacement is written in the form 
\[
\utot(x,t) = \uin(x,t) + u(x,t),
\]
where $u(x,t)$ is the {\it scattered field}.
It is straightforward to show (in the free-space setting) that
$u(x,t)$ must satisfy the following boundary value problem (BVP):
\begin{align}
\pt^2 u - \px^2 u & = 0, \qquad\quad\; x\in\mathbb R\setminus \Gamma, \; t>0,
\label{WE}
\\
u(x,0) = \pt u(x,0) &=  0,\qquad\quad\;  x\in\mathbb{R},
\label{IC}
\\
\jump{u}(x_j,t) & = 0,  \qqquad\!\!\!  j=1,\dots,M,  \; t>0,
\label{cont}
\\
\jump{\px u}(x_j,t) - \scst_j u(x_j,t) &=  g_j(t),  
\qquad j=1,\dots,M, \; t> 0.
\label{kink}
\end{align}
Here the square parenthesis
\[ \jump{s}(x_j,t) \;\equiv\; \lim_{\epsilon \rightarrow 0} \;
s(x_j+\epsilon,t) - s(x_j-\epsilon,t)
\] 
denotes the spatial jump in a function $s$ at position $x_j$ and time $t$.
The time-dependent driving functions are
\beq
g_j(t) := \scst_j \uin(x_j,t) - \jump{\px \uin}{(x_j,t)},
\label{gj}
\eeq
specified by the incident field.
The Cauchy matching conditions at the spring locations
$x_j$ correspond to continuity, and to vertical force balance at
the (massless) point $x_j$. Specifically, the physical string displacement
satisfies
\[ \jump{\px\utot}(x_j,t) = \scst_j \utot(x_j,t), \]
where the left side comes from tension in the string while the right side
is the spring restoring force.
Notice that the limit $\beta_j\to +\infty$ would be a completely reflective Dirichlet condition
allowing no communication between the regions $x<x_j$ and $x>x_j$.
Finite $\beta_j$ gives partial reflection and partial transmission, and is
possibly the simplest nontrivial point scatterer in 1D.%
\footnote{This scatterer can also be realised as the zero-width
limit of a narrow and tall potential function. A related simple point scatterer model is a ``bead on a string'' \cite{vib1983} \cite[Ch.~1]{martinbook}.}
Such scatterers model localized impurities in an otherwise homogeneous medium,
and hence have been studied in the physics and acoustics community
\cite{vib1983,condat86,Disordered2009,Propagation2006}.
Precisely the same ``delta potential'' scatterer model is also
analyzed in quantum mechanics in one or more dimensions \cite{solvableQM,holmer06}.

Our interest in this problem, aside from the physics itself,
is that it will allow us to investigate
integral equation methods based on {\it hyperbolic layer potentials} in
a simple setting. For this, we represent $u(x,t)$ as a 
{\it single layer potential}:
\beq\label{eq:singleLayerSol}
u(x,t) =  {\cal S}[\Gamma,\vec{\sigma}](x,t) 
= \sum_{j = 1}^M\int_0^tG(x - x_j, t - \tau)\sigma_j(\tau)d\tau
=
\frac{1}{2}
\sum_{j = 1}^M\int_0^{t-|x-x_j|} \sigma_j(\tau)d\tau,
\qquad x\in\R, \; t>0,
\eeq
where 
$\vec{\sigma}(t) = \{ \sigma_1(t),\dots,\sigma_M(t) \} $,
each $\sigma_j(t)$ is an unknown {\it density} and
$G(x,t)$ is the 
free-space Green's function for the 1D wave equation:
\beq\label{eq:1DWaveGreens}
G(x,t) = 
\frac{1}{2} H\left(t - \abs{x}\right), 
\eeq
where $H(x)$ as usual denotes the Heaviside step function
\[ H(x) = 
\begin{cases}
1, & x \geq 0 \\
0, & x < 0 .
\end{cases}
\]
With this representation, $u(x,t)$ satisfies \eqref{WE}--\eqref{cont}. Imposing the derivative
jump \eqref{kink} using the jump relation
\beq
[\px u](x_j,t) = -\sigma_j(t)
\label{JR}
\eeq
for the single-layer representation \cite[\S10.6]{Guenther2012},
yields after a little algebra
a system of Volterra integral equations for the unknown densities $\sigma_j(t)$:
\beq\label{eq:IEfree}
-\sigma_j(t) 
- \frac{\scst_j}{2}\sum_{l = 1}^{M} 
\int_0^{t - \abs{x_j - x_l }}\sigma_l(\tau)d\tau
\;=\;
g_j(t), \qquad j = 1, \dots, M, \; t>0.
\eeq
Once the functions $\sigma_j(t)$ are known, the scattered field
can be evaluated at any desired point using
\eqref{eq:singleLayerSol}.

There are several advantages to the integral formulation of
the wave equation. First,
it requires a discretization of the boundary alone (here, the
points $x_j$). 
Second, the order of accuracy of the scheme is 
determined by the order of accuracy of the quadrature in time, 
with no numerical dispersion error. Third, it captures the 
piecewise-smooth structure of the
solution analytically.  

An obvious disadvantage of this approach
is the history-dependence of the layer potentials:
at a given time $t = n \dt$, evaluating the representation \eqref{eq:singleLayerSol}
at a set of $M$ target points---%
or its derivative jumps, as required in the left-hand side of \eqref{eq:IEfree},
would require $\bigO(M^2 n)$ work in a naive implementation. Thus, the 
full simulation would require of the order $\bigO(M^2 N_t^2)$ work.
Moreover, it would require $\bigO(MN_t)$ storage for the history of the densities
at all time steps. 
\begin{remark}
In one space dimension (1D), it is easy to see that $u(x,t)$ satisfies the recurrence
\[
u(x,t+\dt) = u(x,t) +
\sum_{j = 1}^M\int_{t-|x-x_j|}^{t+\dt-|x-x_j|}  \sigma_j(\tau)d\tau. 
\]
This would reduce the cost of evaluating the representation
to $\bigO(M^2)$ work per timestep, for a total work of $\bigO(M^2 N_t)$,
and $\bigO(M/\dt)$ storage. 
Such a simple acceleration is not available in 2D,
due to the weak Huygens' principle,
while in 3D the strong Huygens' principle would yield a similar improvement
without any recurrence (as in, e.g., \cite{Barnett2020}).
In all cases the cost per time step remains quadratic in $M$.
\end{remark}

The main contribution of the present paper 
is a new time-marching method for the wave equation---windowed Fourier projection (WFP)---%
that is independent of spatial dimension,
eliminates the history dependence, requires modest storage,
and whose cost per time step scales quasi-linearly in the number of scatterers $M$.
More precisely, we will make use of a Fourier series with 
$N_F = \bigO(1/\dt)$ modes, in the 1D case, to represent the history contribution
to the representation.
The remaining ``local'' part requires densities only $\bigO(\dt)$ into the past.
The densities and modes evolve in a coupled integro-differential system
that may be efficiently time-marched.
We will show that the total amount of work for our scheme is
$\bigO\bigl( (M+N_F \log N_F) N_t\bigr)$,
assuming that springs are more or less uniformly
distributed and that $\dt \lesssim 1/M$ (these latter conditions ensuring that the
local part has linear cost).

We will first derive a slightly simpler variant of the method
under the assumption that
$\uin(x,t)$ and $u(x,t)$ are spatially periodic on the interval $B$
in the next section.
The evaluation of the history part is discussed in section
\ref{sec:history} and that of the local part in section
\ref{sec:local}, assuming the densities $\sigma_j(t)$
are known. The solution of the integro-differential equation 
to determine the density at each new time step is discussed
in section \ref{sec:inteq} and
the modifications needed for free-space scattering are described
in section \ref{sec:freespace}. 
Section \ref{sec:stab} presents some theoretical analysis of stability and convergence
for simple examples. Section \ref{sec:results} 
presents the results of a variety of numerical experiments.
An appendix contains the proof of the main theorem.

\begin{remark}
Fourier-based methods to overcome the history-dependence
of time-dependent layer potentials were introduced in 
\cite{greengard1990cpam} in the context of the heat equation.
The two key ideas from 
\cite{greengard1990cpam} are:
(a) that high frequency spatial modes are
rapidly damped by the heat kernel, so that a Fourier basis is 
highly efficient, and (b) that the 
coefficients of the spectral representation 
permit a simple evolution formula when marching from one time
step to the next (exploiting the semigroup property of the solution
operator). 

For the wave equation, (b) also holds, but
high frequency modes are {\it not} damped, so that
additional ideas will be needed to make effective use of Fourier analysis.
\end{remark}

\begin{remark}
A number of related approaches 
have been developed over the last few decades.
For the Schr\"{o}dinger equation,
these include the use of
contour deformation methods in Fourier space 
to attenuate high-frequency modes \cite{Kaye2022},
methods that carefully filter outgoing waves at the 
boundary of a finite domain, permitting
the use of the FFT on that domain 
\cite{SOFFER2007,SOFFER2009,SofferStucchio}.
Other methods worth noting that overcome
the quadratic cost of history dependence include
\cite{lubich1985} 
and time-domain versions of the fast multipole 
for the heat and wave equations 
(see \cite{tausch2007jcp} and \cite{fmmbook,Mich2003}, respectively).
\end{remark}

\section{Windowed Fourier projection in the periodic setting} \label{sec:periodic}

To separate the issue of history dependence from that of imposing outgoing 
boundary conditions on a finite computational domain,
we first assume that the incoming and scattered fields are both
periodic on the domain $B = [-\pi,\pi]$.
The scattered field BVP of interest is the wave equation
\eqref{WE} restricted to $x\in B$, with zero initial conditions,
and 
continuity and jump conditions \eqref{cont},\eqref{kink} as before,
plus periodicity imposed by matching Cauchy data:
\beq
u(-\pi,t) = u(\pi,t), \qquad\qquad \px u(-\pi,t) =\px u(\pi,t),  \qquad t> 0.
\label{peri}
\eeq
As before, $\uin(x,t) = f(x-t)$, which controls the driving functions $g_j(t)$ via \eqref{gj}.
However, if $\utot$ is to be periodic, this requires that $f$ be $2\pi$-periodic.%
\footnote{Note that \eqref{WE}--\eqref{kink} with \eqref{peri}
also give a well-posed BVP with arbitrary driving functions $g_j(t)$.}
To solve this BVP, $u$ is represented as a periodic single-layer potential
\beq
\label{SLPp}
u(x,t) 
= \sum_{j = 1}^M\int_0^tG_p(x - x_j, t - \tau)\sigma_j(\tau)d\tau,
\qquad x\in\R, \; t>0,
\eeq
where $G_p(x,t)$ is the periodic 1D wave Green's function on $B$:
\beq
\label{eq:1DperiodicWaveGreens}
G_p(x,t) = \frac{1}{2}\sum_{m= -\infty}^{\infty}H\left(t - \abs{x - 2\pi m}\right). 
\eeq
Rather than \eqref{eq:IEfree}, imposing the jump relations now
yields the system of Volterra integral equations
\beq\label{eq:periodicIE}
-\sigma_j(t) - \frac{\scst_j}{2}\sum_{l = 1}^{M}\sum_{m  = -\infty}^{\infty}\int_0^{t - \abs{x_j - x_l - 2\pi m }}\sigma_l(\tau)d\tau
\;=\; g_j(t),
\qquad j = 1, \dots, M, \; t>0,
\eeq
which has the same history-dependence as before
(in addition to the complexity of having a growing number of image contributions).

We could, however, replace $G_p$
in~\eqref{eq:1DperiodicWaveGreens} by its Fourier series representation:
\beq\label{eq:1DperiodicWaveGreens_Fourier}
G_p(x,t) = \frac{1}{2\pi} \left( t + \sum_{k\neq 0}\frac{\sin kt}{k}
e^{-ikx} \right),
\eeq
and seek the solution in the form
\beq\label{eq:singleLayerSol_Fourier}
u(x,t) = \sum_{k = -\infty}^{\infty}u_k(t)e^{-ikx}.
\eeq
Combining \eqref{SLPp}, \eqref{eq:1DperiodicWaveGreens_Fourier} and \eqref{eq:singleLayerSol_Fourier}
relates the Fourier coefficients $\{ u_k \}$
to the densities $\{ \sigma_j \}$:
\beq\label{eq:fcoef}
u_k(t) =
\int_0^t \frac{\sin k(t - \tau)}{k} {S}_k(\tau) d\tau,
\qquad k\in\Z, \; t>0,
\eeq
where we have introduced the $k$th spatial Fourier coefficient of the density distribution,
\beq
{S}_k(\tau) := \frac{1}{2\pi} \sum_{j = 1}^M \sigma_j(\tau) e^{ikx_j}.
\label{densco}
\eeq
\begin{remark}\label{r:kzero}
  In \eqref{eq:fcoef} the $k=0$ case is taken to mean the limit $\R\ni k\to 0$,
  so that $u_0(t) = \int_0^t (t-\tau) S_0(\tau) d\tau$.
  This limit is implied from now on 
  for the $k=0$ case when $k$ appears in the denominator.
\end{remark}
Note that, since the incident wave is smooth, and the BVP linear,
the response of all the $\sigma_j$, and hence
{all of the ${S}_k$}, are smooth in time.
Furthermore, ${S}_k(\tau)$ ``switches on'' smoothly at $\tau=0$
(i.e., it extends smoothly to the zero function in $\tau\le0$),
because the support of the incident wave at $t=0$ does not include any scatterers.

\begin{figure}[t!]  \label{pkfig} 
  \raisebox{-0.8in}{\includegraphics[width=.67\textwidth]{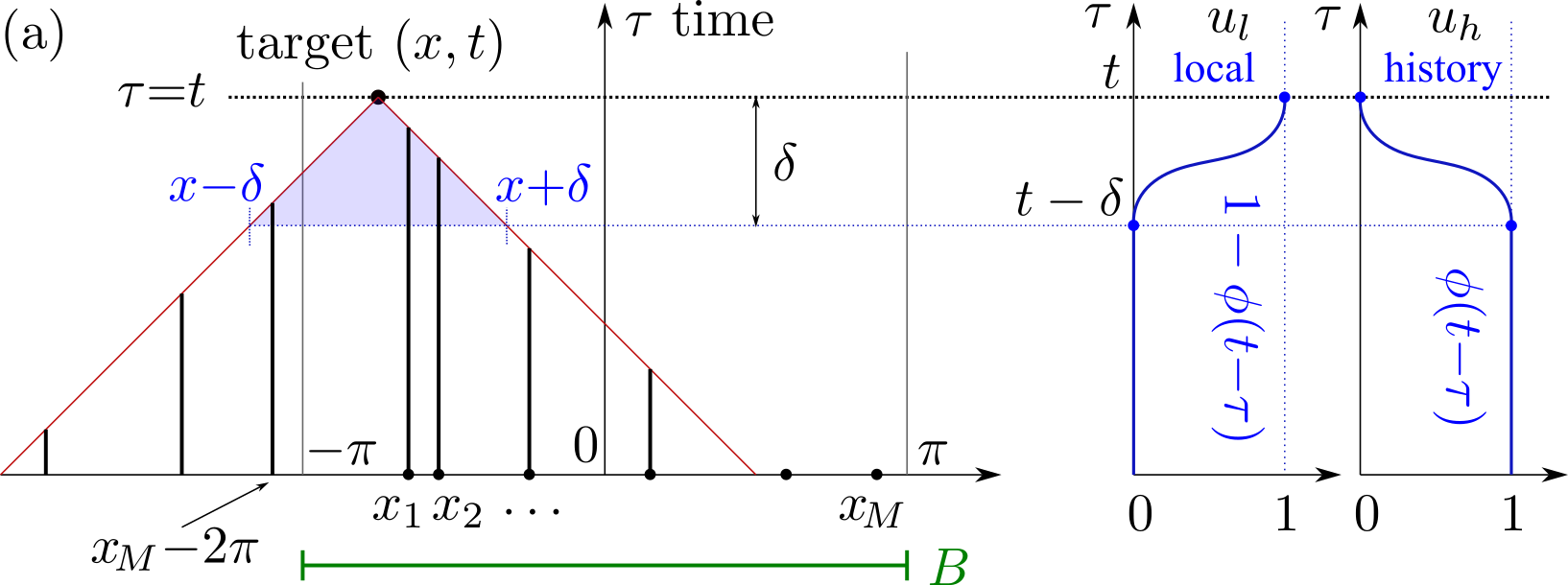}}
\hfill
\begin{minipage}{.3\textwidth}
  (b)
  \\
  \includegraphics[width=\textwidth]{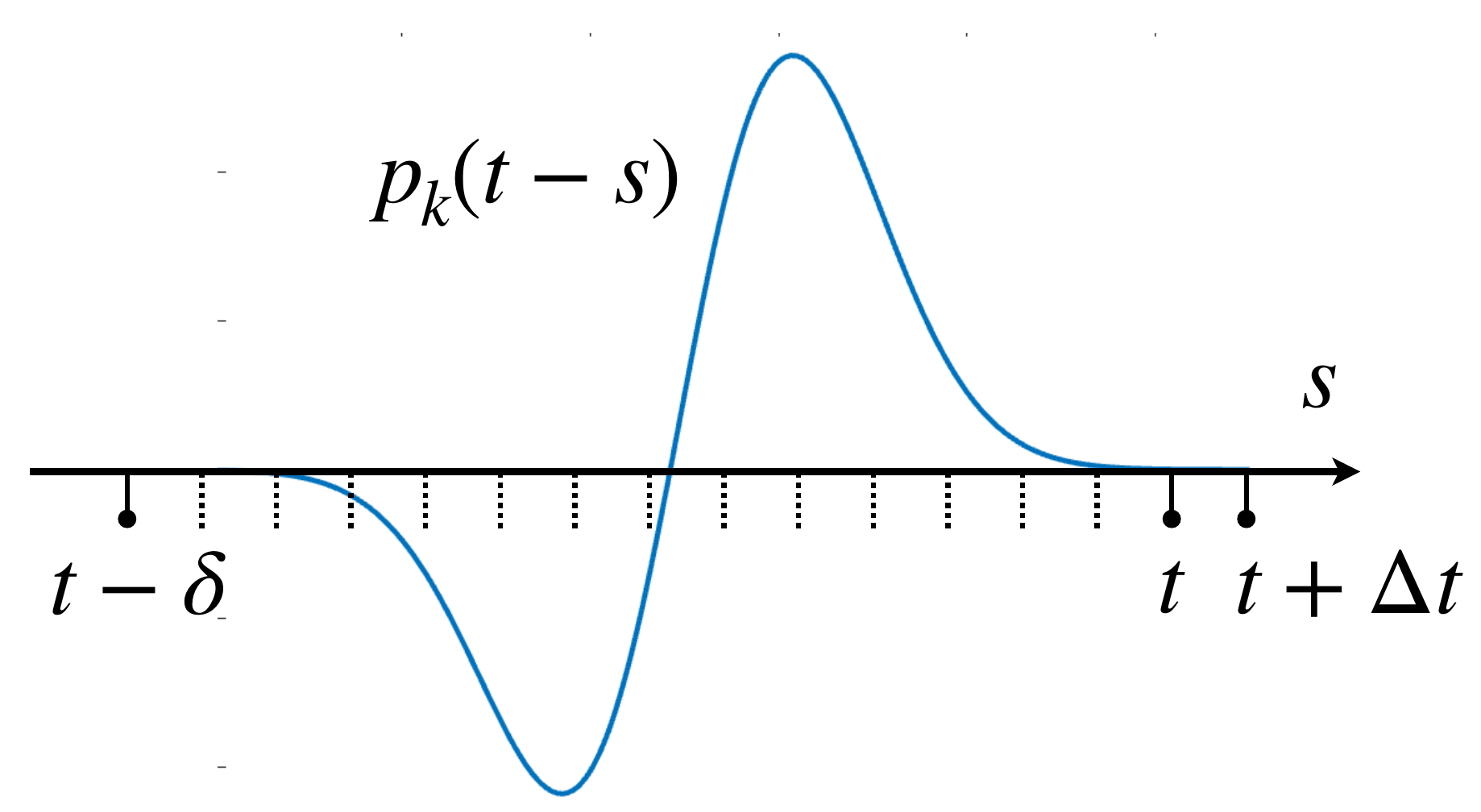}
\end{minipage}
\caption{(a) Space-time diagram of density contributions $\sigma_j(\tau)$
  (vertical lines) for evaluating $u(x,t)$ (left),
  showing the split via the window function $\phi$ into local and history components
  (right). The shaded triangle on the left shows the only spacetime region where densities
  contribute to $u_L$.
  (b) The weight functions
  $p_k$, $q_k$ in \eqref{pkqkcomp}, where $k$ labels Fourier mode,
  vanish at times $s=t-\delta$
and $s=t+\dt$. Thus only known values of 
$\hat\sigma_k(\tau)$ in the recent past $t-\delta \le \tau \le 0$
are involved in computing $h_k(t)$, $g_k(t)$ and,
hence, $\alpha_k(t+\dt)$ (see text for details).
}
\end{figure}

Unfortunately, the Fourier coefficients $u_k(t)$
in \eqref{eq:fcoef} are in general slowly decaying like $\bigO(1/k^2)$, since 
$u(x,t)$ is continuous but only piecewise-smooth
(this low-order decay rate also follows directly from the discontinuous
derivative at the upper limit in \eqref{eq:fcoef}).
By analogy with the heat equation
\cite{greengard1990cpam}, one could try separating the scattered
potential into a local part $u_L$ plus a history part
$u_H$:
\beq
u(x,t) = u_L(x,t) + u_H(x,t),
\eeq
where a sharp transition at some time delay $\delta>0$ into the past is used, so that
\begin{align}
u_L(x,t) &:= \sum_{j = 1}^M\int_{t - \delta}^{t} G_p(x - x_j, t - \tau)
\sigma_j(\tau)d\tau, \label{eq:localSolh}\\
u_H(x,t) &:= \sum_{j = 1}^M\int_0^{t-\delta} 
G_p(x - x_j, t - \tau) \sigma_{j}(\tau) d\tau
.
\label{eq:historySolh}
\end{align}
A simple calculation, however, shows that 
expanding only the history part in a Fourier series does not 
improve its convergence rate.

Here, we introduce the windowed Fourier projection method (WFP),
an alternative decomposition which makes use of 
a smooth partition of unity based on 
a \textit{blending function} $\phi(t)$  with $\phi(t) = 0$
for $t\leq 0$, and $\phi(t) = 1$ for $t\geq \delta$.
We will later set the timescale parameter $\delta$ to be of order several time steps.
Instead of \eqref{eq:localSolh} and \eqref{eq:historySolh}, 
we define the local and history part as
\begin{align}
u_L(x,t) &:= \sum_{j = 1}^M\int_{t - \delta}^{t} G_p(x - x_j, t - \tau)[ 1- \phi(t - \tau)]\sigma_j(\tau)d\tau, \label{eq:localSol}\\
u_H(x,t) &:= \sum_{j = 1}^M\int_0^tG_p(x - x_j, t - \tau)\phi(t - \tau)\sigma_{j}(\tau) d\tau.\label{eq:historySol}
\end{align}
The local part $u_L$ still represents the contribution from the 
densities $\{ \sigma_j(t) \}$ windowed over the most recent time interval (of 
length $\delta$).
In order to evaluate $u_L(x,t)$, we need only include
those springs (or their periodic images) that lie within the interval
$[x-\delta,x+\delta]$, since waves travel at unit speed.
This interval and the windowing are sketched in Fig.~\ref{pkfig}(a).
From now on we assume, for the sake of simplicity, that $\delta < \pi$ so that
no more than one image source can lie in $[x-\delta,x+\delta]$.

\begin{definition} \label{lambdadef} 
For $x \in \R$, we denote by $\dist(x,x_j)$ the distance
from $x$ to the nearest image of the source $x_j$:
\[ \dist(x,x_j) := \min_{m \in \mathbb{Z}} \, |x - (x_j + 2\pi m)|.
\]
We denote by 
$\Nd x$ the set of indices of all points with images lying in
$[x-\delta,x+\delta]$:
\[ \Nd x := \{ j =1,\dots,M \, : \, \dist(x,x_j)  \le \delta \}.
\]
\end{definition}

Using this definition, and recalling \eqref{eq:1DperiodicWaveGreens}, the spatial form of $G_p$,
\beq\label{eq:localSol1}
u_L(x,t) = \frac{1}{2}\sum_{j \in \Nd x}
\int_{t - \delta}^{t-\dist(x,x_j)} 
[ 1- \phi(t - \tau)]\sigma_{j}(\tau)d\tau.
\eeq
An example where $\Nd{x}$ includes a periodically-wrapped source is sketched in Fig.~\ref{pkfig}(a).

For the history part $u_H$ in \eqref{eq:historySol}, we use the Fourier representation
\beq\label{eq:historySol1}
u_H(x,t) = \sum_{k = -\infty}^{\infty}\alpha_k(t)e^{-ikx}.
\eeq
The coefficients $\{\alpha_k\}$ are driven by the
full history of the spatial coefficients ${S}_k$ of the density, recalling \eqref{densco},
except now with time windowing. That is,
\beq
\label{eq:alphak}
\alpha_k(t) :=
\int_0^t \frac{\sin k(t - \tau)}{k} {S}_k(\tau) \phi\left(t - \tau\right)d\tau,
\qquad k\in\Z, \; t>0,
\eeq
which the reader should compare to \eqref{eq:fcoef}.
The key new feature is that ${S}_k(\tau)\phi(t-\tau)$ smoothly
decays to zero at the upper limit $\tau=t$ because of the window function.
Combined with the fact that ${S}_k$ is already known to be smooth and
smoothly vanishing at $\tau=0$,
this implies that for any fixed $t>0$
the integral in \eqref{eq:alphak} vanishes rapidly with respect to
index $|k|\to\infty$: the Fourier series is spectrally convergent.

Making use of the local/history split, 
the coupled system of $M$ Volterra integral equations for 
$\sigma_j$ in \eqref{eq:periodicIE} takes the form
\beq\label{eq:periodicIE_AS}
-\sigma_j(t) 
-\frac{\scst_j}{2}\sum_{l \in \Nd {x_j} }
\int_{t - \delta}^{t - \dist(x_j,x_l)} \!\!\!\!\!
\sigma_{l}(\tau)\left[1 - \phi(t - \tau)\right]d\tau
-
\frac{\scst_j}{2}\sum_{k = -\infty}^{\infty}\alpha_k(t)e^{-ikx_j}
=
g_j(t),  \quad j = 1, \dots, M, \; t>0,
\eeq
where the $\alpha_k$ are given by \eqref{densco} and \eqref{eq:alphak}.
The Volterra system \eqref{eq:periodicIE_AS} with \eqref{densco} and \eqref{eq:alphak}
is equivalent to the original integral equation \eqref{eq:periodicIE} for the
densities alone.

Numerically we will truncate the Fourier series to a maximum index $K$, i.e., we solve
\beq\label{eq:periodicIE_trunc}
-\sigma_j(t) 
-\frac{\scst_j}{2}\sum_{l \in \Nd {x_j} }
\int_{t - \delta}^{t - \dist(x_j,x_l)} \!\!\!\!\!
\sigma_{l}(\tau)\left[1 - \phi(t - \tau)\right]d\tau
-
\frac{\scst_j}{2}\sum_{k = -K}^{K}\alpha_k(t)e^{-ikx_j}
=
g_j(t),  \quad j = 1, \dots, M, \; t>0,
\eeq
where each of the $\alpha_k$ are given by \eqref{densco} and \eqref{eq:alphak} as above.
The spectral convergence rate with respect to $K$ of the truncation error
is controlled by the smoothness of $\phi(t-\tau)$: if the window is less smooth
than ${S}_k(\tau)$, then the $\alpha_k$ will decay excessively slowly as $|k|\to\infty$.
The work in evaluating the Fourier sum in \eqref{eq:periodicIE_trunc} scales like $K$,
so reducing this work suggests a smoother window and a longer $\delta$.
Yet, the cost of the local evaluation (first sum in \eqref{eq:periodicIE_trunc}) grows
linearly with $\delta$, suggesting that one push $\delta$ very small to reduce its work.
These criteria are in conflict.
An efficient scheme thus requires a $\phi$ that is
{\it as smooth as possible for a given blending time $\delta$},
and also chooses $\delta$ to balance the work of local and history parts.

We now formalize this intuition. 
While many windows (blending functions) $\phi$ are possible, in this paper we base
it on the Kaiser--Bessel ``bump'' function
$I_0(b\sqrt{1-t^2})$ on $t\in [-1,1]$, where $b>0$ is a shape parameter, and $I_0$
the modified Bessel function of order zero.
This has the Fourier transform \cite{kaiser,oberhettinger,dftsubmat}
\beq
\int_{-1}^1 I_0(b\sqrt{1-t^2}) e^{i\omega t} dt = 2 \sinc \sqrt{\omega^2 - b^2},
\qquad \omega\in\R,
\label{KB}
\eeq
recalling the definition $\sinc z = \frac{\sin z}{z}$ for $z\neq 0$ and $\sinc 0 = 1$.
For $|\omega|<b$ the right-hand side is exponentially large (nearly $e^b/b$
at $\omega=0$) because
$\sinc iz = \frac{\sinh z}{z}$; in constrast for $|\omega|\ge b$ the right-hand side
has magnitude at most $2$.
Thus the Fourier transform is exponentially (with respect to $b$)
localized%
\footnote{Its degree of localization is close to optimal
in the $L^2$ sense, and has the optimal exponential rate, since the function is
very close to
the prolate spheroidal wavefunction of order zero \cite{slepianI,keranal}.}
to the frequency interval $[-b,b]$.
To create a
blending function $\phi$ transitioning from 0 to 1 over the interval $[0,\delta]$,
we need the antiderivative of a shifted and scaled bump,
\beq
\phi(t) = \int_0^t \phi'(\tau) d\tau,
\quad t\in\R,
\qquad
\mbox{ where }
\quad
\phi'(t) =
\left\{
\begin{array}{ll}
\frac{b}{\delta \sinh b} I_0 \left( b \sqrt{1 - (2t/\delta - 1)^2} \right )
,& \quad 0\le t \le \delta, \\
0, & t< 0 \mbox{ or } t>\delta.
\end{array}
\right.
\label{phi}
\eeq
Although $\phi$ is actually
nonsmooth at $0$ and $\delta$, having derivative jumps of size $\bigO(e^{-b})$,
when one sets
\beq
\label{b}
b = \ln \frac{1}{\epsilon}
\eeq
for some desired tolerance $\epsilon \ll 1$,
then in practice $\phi$ is {\it numerically smooth} to this tolerance.
Using as our definition of the Fourier transform
\beq
\hat{\phi}(\om) := \int_\R \phi(t) e^{i\om t} dt,
\label{FT}
\eeq
we have by rescaling \eqref{KB} (or see \eqref{hatbump})
that $\hat\phi(\om)$ is of size $\bigO(\eps)$
for all $|\om| \ge 2b/\delta$: its $\eps$-support is $[-2b/\delta, 2b/\delta]$.
We refer to $2b/\delta$ as the {\it cutoff frequency} or $\eps$-bandlimit of
the window, recalling the following.

\begin{definition}\label{d:BL}    
  We say a function $f \in L^2(\R)$ is {\it bandlimited} with bandlimit $c>0$ if
  \[
  f(t) = \frac{1}{2\pi}\int_{-c}^c \hat{f}(\om) e^{-i\om t} d\om, \qquad t\in\R.
  \]
  If the difference between left and right sides of the above has $L^2$-norm
  $\bigO(\eps)$ relative to that of $f$, we say that $f$ has $\eps$-bandlimit $c$.
\end{definition}

The principal analytic  contribution of the present paper is the following theorem
which bounds the tail of the Fourier series that has been truncated
in \eqref{eq:periodicIE_trunc}, expressed by the error function
\beq
E_K(t) := \sum_{|k|>K} |\alpha_k(t)|, \qquad t>0.
\label{EK}
\eeq
This controls the {\it backward error} in the following sense.
One may interpret the solution of \eqref{eq:periodicIE_trunc} as the
solution of the true Volterra equation \eqref{eq:periodicIE_AS} with perturbed
data $\tilde g_j(t) = g_j(t) + \frac{\beta_j}{2}\sum_{|k|>K}\alpha_k(t)e^{-ikx_j}$
for each $j=1,\dots,M$.
Then $\beta_j E_K(t)/2$ is an upper bound on the residual $|\tilde g_j(t)- g_j(t)|$.
For wave applications it is common to care about bandlimited excitation $f$
in \eqref{uinc}, and thus by linearity and time-invariance
the solution and the densities must also have the same bandlimit, which we denote
by $K_0$ in the following.

\begin{theorem} \label{t:trunc}  
  Let $\sigma_j$ have bandlimit $K_0$ and $\|\sigma_j\|_\LLR \le C$
  for each $j=1,\dots,M$.
  Let the window $\phi$ be defined by \eqref{phi} with shape parameter \eqref{b}
  and timescale $\delta>0$.
  Let the Fourier series truncation frequency satisfy
  \beq
  K \ge K_0 + \frac{2b}{\delta}.
  \label{Kbnd}
  \eeq
  Then the norm of the error function \eqref{EK} over a finite simulation time $t\in(0,T)$
  has the bound
  \beq
  \|E_K\|_\LLT \le 2\pi^2 TMC \frac{b}{\delta\sinh b}
    =
    \bigO\left(\eps \log \frac{1}{\eps} \right)
    \quad \mbox{ as } \eps\to 0.
    \label{Ebnd}
  \eeq
\end{theorem} 
We defer its proof to \ref{a:pf}.
The big-O statement in \eqref{Ebnd} simply applies \eqref{b}.
We thus expect almost exponential convergence of the residual with respect to $b$
(recalling $\sinh b \sim e^b/2$), as long as the densities are (numerically) bandlimited
and the bandwidth sum condition \eqref{Kbnd} holds.
An interpretation of the latter is that
the windowing (the pointwise multiplication by $\phi$ in \eqref{eq:alphak})
enlarges the density bandlimit $K_0$ by the window $\eps$-bandlimit
$2b/\delta$. The sinusoidal wave propagator term in \eqref{eq:alphak} then converts this
temporal $\eps$-bandlimit to an equal spatial one,
so that Fourier series truncation to this bandlimit induces errors close to $\bigO(\eps)$.
From a physical perspective,
an incident field with bandlimit $K_0$ cannot scatter ``propagating'' spatial frequency
modes of the form $e^{ikx}$ in the actual solution with $|k| > K_0$.
In essence, the singular structure near the springs is resolved by the local
part, while the window function filters out the ``non-propagating''
components from the history part. 

\begin{remark}[Choice of Fourier truncation $K$]\label{r:gamma}  
Choosing an equality in \eqref{Kbnd},
we denote the fraction of the resulting numerical bandlimit $K$ which is ``lost'' to
windowing by $\gamma$, so that $2b/\delta = \gamma K$ while $K=K_0/(1-\gamma)$.
The sharper the window (smaller the timescale $\delta$), the less work
associated with the local part, but the greater the effective bandlimit $K$.
Thus, the selection of the parameter $\delta$ can be used to balance 
the workload between the local and history parts. In the present paper,
we typically fix $\gamma = 1/2$, so that $K = 2K_0$.
In turn, $K_0$ may be chosen as the $\eps$-bandlimit of the incident signal $f$.
\end{remark}

We end the section by mentioning technical limitations in Theorem~\ref{t:trunc}.
1) We were unable to state a hypothesis involving a bandlimit on $f$ or even $g_j$,
and instead made it on $\sigma_j$ (a stronger condition).
We expect (and observe) this to hold to numerical accuracy, but
note that a strict bandlimit is probably not the right condition for a broader analysis.
For instance, bandlimited functions are entire, and thus identically zero
if they vanish for $t<0$: causality and being bandlimited are incompatible.
2) The hypothesis that the densities be in $\LLR$ most likely does not
hold for densities solving the {\it periodic} BVP \eqref{WE}--\eqref{kink} plus \eqref{peri}.
This is because the scattered wave continues to bounce around $B$ for all time.
However, in the more applicable free-space case the densities decay back down to zero
making the condition plausible.
We are unable at this time to weaken this condition to something more local in time,
again because the definition of bandlimited involves all of $\R$.
We leave a refined analysis for future work.

\section{Evaluation of the history part using recurrences} \label{sec:history}

The reason we propose a Fourier representation for the history is 
that, while the Fourier coefficients themselves 
are fully history-dependent,
they can be updated recursively.
A direct computation using \eqref{eq:alphak} shows that
\beq
\alpha_k''(t) + k^2 \alpha_k(t) = F_k(t),
\qquad \mbox{ with driving }
\quad
F_k(t):=\int_{t-\delta}^t
\Psi_k(t-\tau) {S}_k(\tau) d\tau,
\label{Fk}
\eeq
given in terms of the blending ``influence kernel''
\beq
\Psi_k(\tau) := 2\cos k\tau \, \phi'(\tau) +\frac{\sin k \tau }{k}\phi''(\tau),
\qquad \tau \in \R,
\label{psidef}
\eeq
whose support is $[0,\delta]$.
Thus the $k$th coefficient obeys a 2nd-order linear ODE driven by $F_k$, a weighted
integral of ${S}_k$ over the most recent time interval of length $\delta$.
It is useful to rewrite \eqref{Fk} as the linear 1st-order system
\beq
\partial_t
\begin{bmatrix}\al_k(t) \\ \al_k'(t) \end{bmatrix}
= \begin{bmatrix}0 & 1 \\ -k^2 & 0 \end{bmatrix}
\begin{bmatrix}\al_k(t) \\ \al_k'(t) \end{bmatrix} +
\begin{bmatrix}0 \\ F_k(t) \end{bmatrix},
\qquad t>0.
\label{FOS}
\eeq
For each $k\in\Z$, the matrix which propagates the homogeneous system over a time $t$
is easily seen to be
\beq
\begin{bmatrix}\cos k t & \frac{\sin kt}{k} \\ -k \sin kt & \cos kt \end{bmatrix}.
\eeq
The Duhamel principle then proves the following.
\begin{lem}  
Let $\Delta t$ denote a time step. Then
the Fourier coefficients $\alpha_k(t)$ satisfy the evolution formulae
\begin{equation}\label{eq:alphakEvolutionFormulas}
\begin{split}
\alpha_k(t + \dt) &= \alpha_k(t)\cos k\dt + \alpha_k'(t)\frac{\sin k\dt}{k} + h_k(t), \\
\alpha_k'(t + \dt) &= -k\alpha_k(t)\sin k\dt + \alpha_k'(t)\cos k\dt + g_k(t),  
\end{split}
\end{equation}
where, recalling the definition of $F_k$ in \eqref{Fk} and \eqref{psidef},
\beq \label{eq:hkgk}
h_k(t) := \int_t^{t + \dt}\frac{\sin k(t + \dt - \tau)}{k}F_k(\tau)d\tau,
\qquad
g_k(t) := \int_t^{t + \dt}\cos k(t + \dt - \tau)F_k(\tau)d\tau.
\eeq
\end{lem}

The preceding lemma yields exact
one-step recurrences for $\alpha_k$, and $\alpha'_k$. 
The driving updates, $h_k(t)$ and $g_k(t)$,
however, need to be computed numerically.
From \eqref{Fk} we have
\beq
\begin{aligned}
h_k(t) &= \int_t^{t + \dt}\frac{\sin k(t + \dt - \tau)}{k}
\left\{\int_{\tau-\delta}^\tau
\Psi_k(\tau-s) {S}_k(s) ds \right\} \, d\tau,
\\
g_k(t) &= \int_t^{t + \dt}\cos k(t + \dt - \tau)
\left\{\int_{\tau-\delta}^\tau 
\Psi_k(\tau-s) {S}_k(s) ds \right\} \, d\tau
.
\end{aligned}
\eeq
Because $\Psi_k$ vanishes outside $[0,\delta]$, one
can change the limits of the inner integrals to $[t-\delta,t+\dt]$,
then swap the order of integration to yield
\beq \label{hkgkcomp}
h_k(t) = \int_{t-\delta}^{t+\Delta t}
p_k(t-s)  {S}_k(s) ds,
\qquad
g_k(t) = \int_{t-\delta}^{t+\Delta t}
q_k(t-s) {S}_k(s) ds,
\eeq
where the ``one timestep integrated'' influence kernels are 
\beq \label{pkqkcomp}
\begin{aligned}
p_k(t-s) &:=
\int_0^{\Delta t} \frac{\sin k (\Delta t - \mu)}{k} \Psi_k(t-s+\mu) \, d\mu,
\\
q_k(t-s) &:= 
\int_0^{\Delta t} \cos k (\Delta t - \mu)\, \Psi_k(t-s+\mu) \, d\mu.
\end{aligned}
\eeq
Recalling from \eqref{phi} that $\Psi_k$
is smooth on $\R$ apart from discontinuities of size $\bigO(\epsilon)$
at either end of its support $[0,\delta]$,
we see that $p_k$ and $q_k$
vanish at the endpoints $t-s=-\dt$ and $t-s=\delta$,
and are smooth apart from derivative discontinuities of size $\bigO(\epsilon)$
(see Fig.~\ref{pkfig}(b)).
As a result, the trapezoidal rule  \cite{PTRtref} is {\it spectrally accurate}
(up to $\bigO(\epsilon)$)
for approximating the integrals in \eqref{hkgkcomp},
and we thus use the available uniform time grid,
with time step $\Delta t$, as its quadrature nodes.
Specifically, letting $n$ be the index of the current time step, so that $t = n\dt$,
we employ the discrete formulae
\beq  \label{hkcompdisc}
\begin{aligned}
h_k(n\dt) &\approx  \dt \sum_{m = 0}^{W-1}
p_k(m\dt) \, {S}_k \left( (n-m)\dt\right), \\
g_k(n\dt) &\approx  \dt \sum_{m = 0}^{W-1} 
q_k(m\dt) \, {S}_k\left( (n-m)\dt\right),
\end{aligned}
\eeq
where we chose
\beq
\delta = W \dt
\label{W}
\eeq
for some integer $W$ setting the number of time steps in the window width.
These explicit formulae, involving only known data on the grid,
are inserted into \eqref{eq:alphakEvolutionFormulas} to
time-step the $\al_k$ and $\al'_k$ coefficients.
In the above, the needed terms $p_k(m\Delta t)$ and
$q_k(m\Delta t)$ for $m=0,\dots,W-1$ are easily approximated
to negligible error by applying
a $N_g$-node Gauss--Legendre quadrature rule to the integrals \eqref{pkqkcomp}.

\subsection{Parameter choices and a fast algorithm for the history part}

The windowed Fourier projection and the timestepping of its history part described
so far is a spectral method.
Recall that $K$ is the highest frequency component kept in the
truncated Fourier series, so that there are a total of $N_F = 2K+1$ modes.
A grid with time step $\dt$ cannot discretize signals beyond a bandlimit
of $\pi/\dt$ (i.e., the signal with alternating signs on the grid)
without aliasing; this is the {\it Nyquist frequency}
\cite[\S4.6]{dahlquistv1}.
This gives an upper bound on $K$ (the same upper bound results
by requiring that the trapezoid quadrature be accurate out to the
maximum bandlimit $2K$ of the integrands in \eqref{hkgkcomp}).
Making $K$ larger than this is waste.
Thus, we connect the spatial and temporal discretizations via
\beq
K = \frac{\pi}{\dt},
\label{Kdt}
\eeq
based on whatever $\dt$ is chosen for accuracy and efficiency reasons.
We now see that $W = \bigO(1)$, the window timescale in time steps, as follows.
Remark~\ref{r:gamma} stated that the window $\eps$-bandlimit is $2b/\delta = \gamma K$,
and combining with \eqref{Kdt} and \eqref{W} gives
\[
W = \frac{2b}{\pi \gamma} = \frac{2}{\pi\gamma} \log \frac{1}{\eps}.
\]
For instance, $\eps=10^{-12}$ and $\gamma=1/2$ gives $W \approx 35$, a setting that
we typically use in later experiments.

We now turn to computational cost. Computing the values 
$p_k(m\Delta t)$, $q_k(m\Delta t)$  for
$m = 1,\dots,W$ requires $\bigO(KWN_g)$ work, where $N_g$
denotes the number of Gauss--Legendre nodes used for quadrature of \eqref{pkqkcomp}.
Since $W$ and $N_g$ are of order 1, this precomputation
is effectively $\bigO(K)$ work and is negligible.

At the $n$th time step, computing
$h_k(n\dt)$ and $g_k(n\dt)$ via \eqref{hkcompdisc} involves two steps.
First we must compute ${S}_k(n\Delta t)$
according to \eqref{densco}, repeated here:
\[
{S}_k(n\dt) := \frac{1}{2\pi} \sum_{j = 1}^M \sigma_j(n\dt)e^{ikx_j},
\qquad k=-K,-K+1,\dots, K.
\]
A naive calculation would require $\bigO(MK)$ work.
However, the form of the exponential sum is that of a nonuniform FFT
of type 1 \cite{nufft2,nufft3,nufft6}, so the set of values $\{ {S}_k(n\dt)\}_{|k|\le K}$ 
can be approximately computed together to relative error $\eps$
in $\bigO\left(M\log(1/\eps)+ K \log K\right)$ time.
These values for the previous $W$ time steps must be stored
to compute \eqref{hkcompdisc}, needing $\bigO(WK)$ storage.
Then, applying \eqref{eq:alphakEvolutionFormulas} costs only $\bigO(K)$.
Thus, obtaining the $\alpha_k$ coefficients of the history part
\eqref{eq:historySol1} at each time step requires $\bigO(M + K \log K)$ work.

The other history-related task each time step
is to evaluate the history representation back at the scatterers, 
for the $\alpha$-dependent term in \eqref{eq:periodicIE_trunc}.
This takes the form of a nonuniform FFT of type 2;
it is the adjoint of the type 1 and requires the same computational effort.
If other solution evaluation points
are needed, they can be appended to the list of scatterers in this step,
with a linear cost.
In summary, handling the history part for the Volterra
equation \eqref{eq:periodicIE_trunc} requires $\bigO(M + K \log K)$ work per time step.

\begin{figure}[t]  \label{lightconefig}  
\centering
\includegraphics[width=\textwidth]{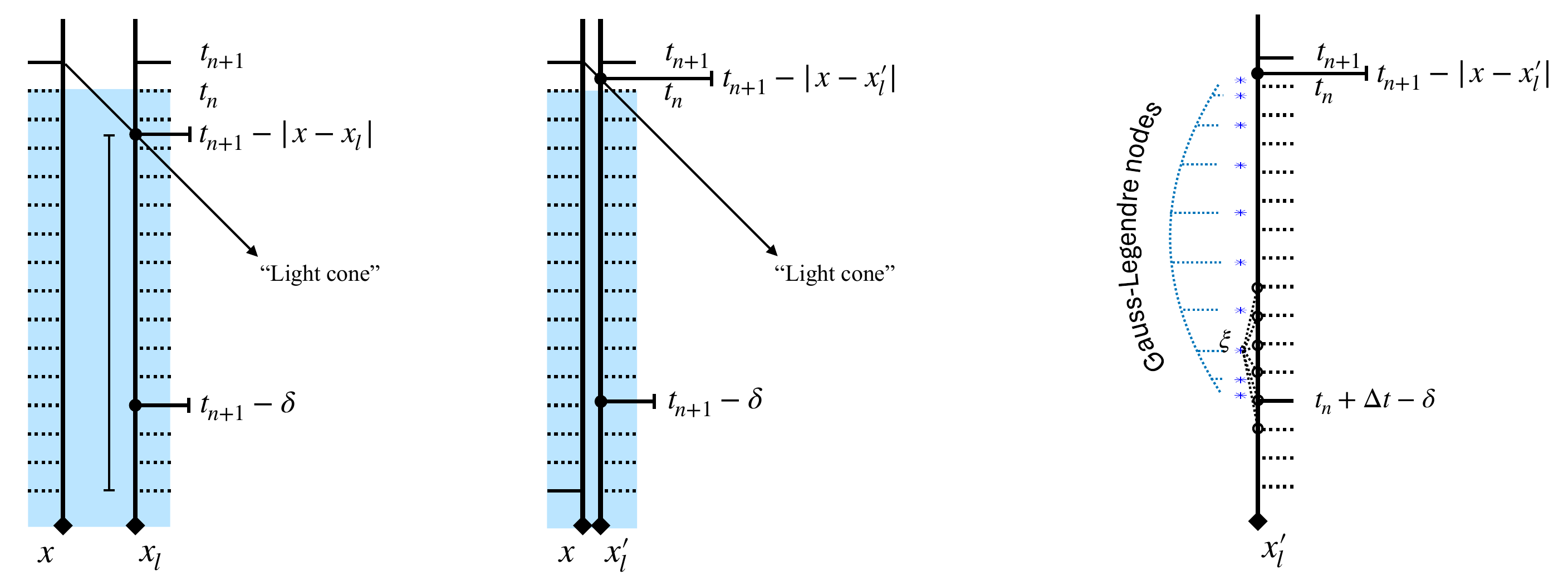}
\caption{Left side: source point $x_l$ has
  ``intermediate'' distance from $x$, i.e., between $\dt$ and $\delta$,
so the characteristic (light cone) emanating from spacetime point $(x,t_{n+1})$
hits $x_l$ before $t_n$. Thus, the integrand
over the interval $[t_{n+1}-\delta,t_{n+1} - |x-x_l|]$ 
can be computed from known
values (indicated by the shaded region) of the density $\sigma_l(t)$.
Middle: the source point $x_{l'}$ is ``nearby'' (closer than $\dt$),
so that the light cone intersects $x_{l'}$ {\it after} the time $t_n$ (outside the
shaded region), and interpolation requires the unknown density $\sigma_l(t_{n+1})$.
Right: quadrature nodes $\xi_i$ for the integration domain
$[t_{n+1}-\delta, t_{n+1}-\dist(x,x_l)]$, in the nearby case.
The small fan of dotted lines indicates degree-$(p-1)$ local interpolation to a single node
$\xi$ from the $p$ nearest regular grid values.
} 
\end{figure}

\section{Evaluation of the local part, time marching, and free space boundary conditions}
\label{s:method}

In this section we describe the remaining ingredients needed to
efficiently time march a fully discrete version of the Volterra integral equation system
\eqref{eq:periodicIE_trunc}, and to change the boundary conditions from periodic
to free-space.

\subsection{Approximation of the local part via $p$th-order density interpolation}
\label{sec:local}

Here we tackle the local part of the single layer
potential, namely the 2nd term in \eqref{eq:periodicIE_trunc}, 
needed to advance from time $t_n = n\dt$ to $t_{n+1} = (n+1)\dt$.
We denote the target point by $x$; this could be the $j$th scatterer $x_j$,
or another point at which the scattered field $u(x,t_{n+1})$ is desired.
We distinguish ``nearby'' sources with $\dist(x,x_j) < \Delta t$
from the remaining more numerous
``intermediate'' distance sources with $\Delta t \le \dist(x,x_j) < \delta$.
Recalling Definition~\ref{lambdadef}, the nearby indices
are $\Ndt x$, while we denote the intermediate indices by
$\mathcal{I}_{\delta}(x) := \Nd x \backslash \Ndt x$.
From \eqref{eq:localSol1}, using $l$ as the source index, we have
\beq\label{eq:localSol2}
u_L(x,t_{n+1}) =
\!\!
\sum_{l \in \Ndt x}
\int_{t_{n+1} - \delta}^{t_{n+1}-\dist(x,x_l)} 
\!\!\![1- \phi(t_{n+1} - \tau)]\frac{\sigma_l(\tau)}{2}d\tau
+
\!\!
    \sum_{l \in \mathcal{C}_{\delta}(x)}
\int_{t_{n+1} - \delta}^{t_{n+1}-\dist(x,x_l)}
\!\!\!
[1- \phi(t_{n+1} - \tau)]\frac{\sigma_l(\tau)}{2}d\tau.
\eeq
The reason for distinguishing the two sets of source points is that
computing the integrals for $l \in \Ndt x$ will involve 
the unknown $\sigma_{l}(t_{n+1})$
(see Fig.~\ref{lightconefig} middle and right panels),
while the integrals for $l \in \mathcal{C}_{\delta}(x)$ use only known values
of the density $\sigma_l$ at time $t_n$ or earlier
(Fig.~\ref{lightconefig} left panel).

Each integral in the sums in \eqref{eq:localSol2} is approximated by
a single Gauss--Legendre quadrature rule over their domain,
with the off-grid density values at these quadrature nodes approximated
via local degree-$(p-1)$ polynomial interpolation from the
$p$ nearest known regular grid values.
Specifically, consider the term involving
either a nearby or intermediate source $x_l$, and
let $\{\xi_i\}_{i=1}^{N_L}$ be the Gauss--Legendre nodes and $\{w_i\}_{i=1}^{N_L}$
the corresponding weights for the interval $\tau \in [t_{n+1} - \delta, t_{n+1}-\dist(x,x_l)]$.
The interpolant evaluated at one of these nodes $\xi$
depends linearly on the grid values $\sigma_l^{n} := \sigma_l(t_n)$, as follows:
\beq
\tilde\sigma_l(\xi) = \sum_{m\in{\mathcal I}(\xi)} v_m(\xi) \sigma_l^{n+1-m},
\label{interp}
\eeq
where ${\mathcal I}(\xi)$ is the set of indices of the $p$ nearest available points to $\xi$
from the grid $\{t_{n+1-m}\}_{m=\mmin}^\mmax$.
See Fig.~\ref{lightconefig}, right panel.
Here  $\mmax = W + \lceil p/2 \rceil$, to allow for interpolation to target
times the furthest in the past, $t_{n+1-W}$.
The lower allowable index is $\mmin = 0$ (accessing $t_{n+1}$) if $x_l$ is a nearby source, or
$\mmin = 1$ if an intermediate distance source.
The interpolation weights $v_m$ are computed by standard barycentric Lagrange interpolation
\cite[\S4.2.3]{dahlquistv1}.
We expect an accuracy $\bigO((\dt)^p)$ for this $p$-node polynomial interpolant.

Each integral term in \eqref{eq:localSol2} is then approximated as
\[
\int_{t_{n+1} - \delta}^{t_{n+1}-\dist(x,x_l)}
\!\!\! [ 1- \phi(t_{n+1} - \tau)]\sigma_{l}(\tau)d\tau
\approx
\sum_{i=1}^{N_L} w_i [ 1- \phi(t_{n+1} - \xi_i)]\tilde\sigma_{l}(\xi_i)
\approx
\sum_\mmin^\mmax Q_m(\dist(x,x_l))\sigma_l^{n+1-m}
.
\]
The formula for the weights is, from \eqref{interp},
\[
Q_m(\dist(x,x_l)) = \sum_{i : m\in{\mathcal I}(\xi_i)} w_i [ 1- \phi(t_{n+1} - \xi_i)] v_m(\xi_i),
\qquad m=\mmin,\mmin+1,\dots,\mmax,
\]
recalling that the nodes $\xi_i$ and weights $w_i$ themselves depend on $\dist(x,x_l)$.
If $\dist(x,x_l)>\delta$, we define all $Q_m$ as zero.

The above holds for local evaluation at an arbitrary target $x\in B$.
For the specific case of $x=x_j$ we will abbreviate
$Q_{jlm}^{\cal N}:=Q_m(\dist(x_j,x_l))$ when $x_l$ is nearby, while
$Q_{jlm}^{\cal C}:=Q_m(\dist(x_j,x_l))$ when $x_l$ is at intermediate distance.
The complete formulae for $Q_m$ incorporating barycentric weights
are cumbersome and available in the accompanying software \cite{springsoftware}.

\subsection{The fast method for time marching}
\label{sec:inteq}

Let us abbreviate $g_j^{n} := g_j(t_n)$, and $\alpha_k^n = \alpha_k(t_n)$, $|k| < K$,
for the time steps $t_n = n\dt$, $n = 0, 1, \dots, \Nt$. 
A fully discrete version of the
integral equation \eqref{eq:periodicIE_trunc} is,
setting $t=t_{n+1}$ and applying the above local evaluation quadratures,
\beq\label{eq:periodicIE_disc}
-\sigma_j^{n+1}
-\frac{\scst_j}{2}
\sum_{l \in \Ndt {x_j} }
\sum_{m=0}^\mmax Q^{\cal N}_{jlm} \sigma_l^{n+1-m} 
-\frac{\scst_j}{2}
\sum_{l \in \mathcal{C}_{\delta}(x_j)}
\sum_{m=1}^\mmax Q^{\cal C}_{jlm} \sigma_l^{n+1-m}
-
\frac{\scst_j}{2}\sum_{k = -K}^{K}\alpha_k^{n+1} e^{-ikx_j}
=
g_j^{n+1},
\eeq
for $j=1,\dots,M$.
Notice that the nearby (first) sum involves $\sigma_l^{n+1}$, and thus is
an ``implicit'' formula,
while the intermediate (second) sum plus the history part
(since the $\alpha_l^{n+1}$ update only needs densities at $t^{n}$ or earlier),
are ``explicit'' terms.
This yields an implicit marching scheme that requires the solution of a 
sparse, banded linear system at each time step, with the number of 
non-zeros in the $j$th row determined by 
the number of nearby points, namely $|\Ndt {x_j}|$.
(If no scatterers are within $\Delta t$ of each other, 
the system of equations is diagonal.) 

Defining at the $n$th time step the vectors of densities
$\bm \sigma^n := [\sigma_1^n, \dots, \sigma_M^n]$,
right-hand sides $\bm g^n := [g_1^n, \dots, g_{M}^n]$, and 
history coefficients $\bm \alpha^n := [\alpha_{-K}^n, \dots, \alpha_{K}^n]$,
it is easy to see that
\eqref{eq:periodicIE_disc} can be written as a linear system of the form
\beq\label{eq:linearSystem}
-(I + S)\bm{\sigma}^{n+1} = \bm g^{n+1} + \sum_{m = 1}^\mmax C^{(m)}\bm\sigma^{n+1-m} + E\bm{\alpha}^{n+1}
,
\eeq
where $\bm\al^{n+1}$ is obtained via the history update equations
\eqref{eq:alphakEvolutionFormulas} and \eqref{hkcompdisc}, implemented via the type 1 NUFFT
as in Section~\ref{sec:history};
these only involve $\bm\sigma^{n+1-m}$ for $m=1,\dots,W$, and $\bm\al^{n}$ and ${\bm\al^{n}}'$.
Thus all of the terms on the right-hand side of \eqref{eq:linearSystem} are known.
The system matrix has the following structure:
$I$ is the $M\times M$ identity matrix, while $S$ is a sparse matrix 
wth non-zero entries $S_{jl} = \frac{\scst_j}{2} Q^{\cal N}_{jl0}$
whenever $\dist(x_j,x_l)< \Delta t$, for $j,l = 1,\dots,M$.
Each matrix $C^{(m)} \in \mathbb{R}^{M \times M}$ is sparse
wth non-zero entries
\[
C^{(m)}_{jl} =
\left\{\begin{array}{ll}
(\beta_j/2) Q^{\cal N}_{jlm}, & \dist(x_j,x_l)<\dt, \\
(\beta_j/2) Q^{\cal C}_{jlm}, & \dt\le \dist(x_j,x_l)<\delta,
\end{array}\right.
\qquad \mbox{ for } \; j,l \in 1,\dots, M.
\]
$E$ represents the discrete nonuniform Fourier transform matrix of
size $M\times (2K+1)$ with entries
\beq
E_{jk} = \frac{\beta_{j}}{2} e^{-ikx_j},
\eeq
and can be applied rapidly using the type 2 NUFFT, as discussed at the end of
Section~\ref{sec:history}.

We expect the error for one time step to scale like the interpolation
order $p$ plus one, since the local integral has an extra factor $\dt$.
The history part we expect to be spectrally accurate, with parameters
chosen for a fixed tolerance $\eps$.
Since there are $N_t = \bigO(1/\dt)$ time steps to reach a fixed final time $T$,
the overall error we expect to be at least $\bigO((\dt)^p + \eps)$.

\subsection{Computational complexity of the full method}
\label{s:compfull}

From now on let $\ntyp:= (1/M)\sum_{j=1}^M |\Nd {x_j}|$ denote the
mean number of neighbors involved in the local part evaluation of a scatterer.
If the $\dt$ needed to achieve accuracy for the given incident signal $f$ is too large,
we recommend decreasing $\dt$ for efficiency reasons so that $\ntyp$ is $\bigO(1)$.
In practice we find that $\ntyp$ of order $10^2$ balances history and local costs.
Precomputing all local quadrature weights $Q^{\cal N}_{jlm}$ and $Q^{\cal C}_{jlm}$
is a one-time cost of $\bigO(N_L p \ntyp M) = \bigO(M)$, with a large prefactor
(however we will see that its cost is less than $10^2$ time-steps).
By a convergence study we found that $N_L = W$ is sufficient. 
The cost of evaluating the local part $u_L(x_j,t_{n+1})$ at a single target is $\bigO(W \ntyp)$,
thus its total cost per time step is $\bigO(M)$.
In practice we reduce its prefactor by applying the stack of $C^{(m)}$ matrices via a single
sparse matrix-vector multiply.
The storage cost for the densities over the needed recent
past is $M\mmax = M(W+ \lceil \frac{p}{2} \rceil)$ complex numbers.
We also store the recent $W$ vectors of $\{ {S}_k\}_{k=-K}^K$ to avoid them having
to be recomputed. This means that there is a one
type 1 and one type 2 NUFFT, a cost $\bigO(M + K \log K)$, each time step.

As for the integral equation itself, to obtain
$\bm{\sigma}^{n+1}$ requires solving \eqref{eq:linearSystem} each time step.
This requires $\bigO(M)$ work per time step since $S$ is banded with a bandwidth
of $\max_j \Ndt {x_j}$, which is $\bigO(\ntyp)$ if the
springs are more or less uniformly distributed.
We accelerate the solve by sparse LU factorization of $I+S$ once and for all.

Note that the condition $\ntyp=\bigO(1)$ means that $\dt = \bigO(1/M)$,
and recalling \eqref{Kdt}, $K=\bigO(M)$.
Thus the overall cost per time step may be summarized as $\bigO(M + K \log K) = \bigO(M \log M)$,
giving a total solution cost over a fixed interval $(0,T)$ of
$\bigO(N_t M \log M) = \bigO(M^2 \log M)$.
This contrasts with $\bigO(N_t^2 M^2)$ for naive time-stepping of the Volterra equation.

\subsection{Imposing free-space boundary conditions via outgoing windowing}
\label{sec:freespace}

Imposing radiation boundary conditions (RBCs) for the scattered solution of the wave equation
in the 1D box $B$ is
easy compared to the higher dimensional case,
because they are {\it local} in space and time, namely $\px u \pm \pt u = 0$
at $x=\pm \pi$ respectively.
(For some of the recent literature on the higher-dimensional case
see \cite{hagstromnrbc,tsynkov01,SOFFER2009,SofferStucchio}.)
In the spirit of the presented WFP method, we present a simple---but we believe new---spectrally accurate
method which acts on the history part.
At certain time steps, a fixed matrix
$P \in \C^{2N_F \times 2N_F}$ modifies the coefficient vectors
$\bm\al := \{\al_k\}_{k=-K}^K$ and $\bm\al' := \{\al'_k\}_{k=-K}^K$ via
\[
\begin{bmatrix}\bm\al\\ \bm\al' \end{bmatrix}
\mapsfrom
P
\begin{bmatrix}\bm\al\\ \bm\al' \end{bmatrix}.
\]
We will see that this must be performed every few time steps.

The action of $P$ is best described as the composition of three steps.
Step 1 takes the FFT of the vectors $\bm\al$ and $\bm\al'$ to get
values on a spatial grid.
It is slightly neater (and incurs negligible error) to first discard the highest frequency
$k=K$, so that the FFTs have length $n=2K$.
The FFT of the (FFT-shifted) vector
$\bm\al$ then outputs the grid of values
$u_H(0), u_H(h),\dots,u_H((K-1)h), u_H(-Kh),u_H((-K+1)h) \dots,u_H(-h)$,
where the grid spacing is $h = \dt = 2\pi/n$, recalling \eqref{Kdt}.
We will drop the $H$ (history) subscript.
The FFT of $\bm\al'$ generates a similar grid but for the function $v := \pt u$.

For step 2
we now interpret the window $\phi$ given by \eqref{phi}, with the same shape parameter $b$
chosen for tolerance $\eps$ as in \eqref{b}, as a
spatial rather than temporal blending function. The blending occurs
from 0 up to 1 over $W$ gridpoints.
Consider imposing a right-going condition at the right endpoint at $x=\pi$.
The (zero-indexed) index set $I = \{K-2W,K-2W+1,\dots,K-1\}$ in the $u$ grid from step 1
corresponds to the interval $R = [\pi-2\delta,\pi)$ where modification will occur.
At $x = ih$ for each $i\in I$ we compute ``rolled-off'' grid values
\begin{align}
  \tilde u(x) &= \phi(\pi-\delta-x) u(x),
  \\
  \tilde v(x) &= \phi'(\pi-\delta-x) u(x) + \phi(\pi-\delta-x) v(x).
  \label{vtilde}
\end{align}
The window $\phi(\pi-\delta-x)$ starts at 1 at the left end of the region $R$, drops to 0 by
the middle of $R$, and stays zero until its right end.
Thus $\tilde u$ is a version of $u$ rolled off to zero over the distance $\delta$.
If $u$ and $v=\pt u$ already describe a right-going wave
in $R$, then \eqref{vtilde} sets $\tilde v$ to have the
right-going condition to match $\tilde u$:
\[
\tilde v(x) = \pt \tilde u(x) = -\px \tilde u(x)
= -\px (\phi(\pi-\delta-x) u(x)) =
\phi'(\pi-\delta-x) u(x) - \phi(\pi-\delta-x) \px u(x),
 \]
and finally the right-going hypothesis $\px u = -\pt u$ justifies \eqref{vtilde}.
The $u$ and $v$ values at the grid point indices $I$ are now replaced with $\tilde u$ and
$\tilde v$.
A reflection of Step 2 is similarly applied on the left end, for which the index set is
$I = \{K,K+1,\dots,K+2W-1\}$.

Step 3 applies inverse FFTs to the resulting $u$ and $v$ grids to give the
new coefficient vectors $\bm\al$ and $\bm\al'$, appending by zero for the $k=K$ entry.
This completes the action of $P$.
The MATLAB implementation is about 10 lines of code (see {\tt alphaRBCproj.m} in \cite{springsoftware}).

We now explain the requirement on the period $\dt_{\text{proj}}$ between
such RBC ``projections.''\footnote{We use quotes here since $P^2 \neq P$,
associated with $\phi$ having a continuous transition from $0$ to $1$.}
The new coefficients represent a field unchanged in $|x|<\pi-2\delta$, but
rolled off to zero and outgoing in $|x| \in [\pi-2\delta, \pi]$.
One can wait up to $2\delta$
before the support of the rolled off wave wraps around the periodic domain and
hits $|x| =\pi-\delta$, where it would start to break the hypothesis that
$(u,v)$ represent outgoing waves in $|x| \in [\pi-2\delta,\pi-\delta]$.
Thus $\dt_{\text{proj}} \le 2\delta = 2W\dt$.
However, repeating before $\delta$ would cause repeated window
  multiplication on the same piece of traveling wave,
unnecessarily growing its spectral bandwidth and perhaps causing aliasing.
In summary, any $\dt_{\text{proj}}$ between $W$ and $2W$ time steps is valid.
An increase in bandwidth by $\gamma K$ is expected (Remark~\ref{r:gamma});
however, this only affects waves that will shortly be zeroed out, and
empirically we find errors $\bigO(\eps)$ matching the window tolerance.

To ensure that no local evaluation affects either rolled-off region, the
scatterers must all reside in a computational domain $\Omega := [-\pi+3\delta, \pi-3\delta]$,
only slightly smaller than the periodic box $B = [-\pi,\pi]$.

The cost of this RBC is four FFTs of length $N_F$, performed every
$\bigO(W)$ time steps, negligible compared to the integral equation cost.
The complexity of the free-space method is thus the same as that of the spatially periodic
method.

\section{Stability and convergence for model problems}
\label{sec:stab}
Here we prove some results for $p$th-order interpolatory
marching as in Sections \ref{sec:local}--\ref{sec:inteq},
applied to
\eqref{eq:IEfree}, the original system of Volterra integral equations involving the free-space
Green's function.
This is the relevant question, because the WFP (local/history split), RBC,
and Gauss--Legendre quadratures together
approximate \eqref{eq:IEfree} to spectral accuracy
(their errors that can easily be pushed down to close to machine precision).

We study discretizations coming from 
1st- or 2nd-order density interpolants, and the case of $M=1$ and $M=2$
springs (scatterers); we believe that this covers the main issues arising with 
nearby or well-separated springs.
It also serves to illustrate the time stepping schemes in simple cases.
Since we seek stability results that do not allow exponential growth in the prefactor,
we employ direct methods as in \cite[\S7.5]{linzbook}.
The stability proofs already involve subtleties, so we defer more generality to future work.
We end with a 2nd-order convergence proof for two springs in the style of
\cite[Thm~7.2]{linzbook},
which in contrast will involve an exponentially growing prefactor.

\subsection{Stability for one spring ($M=1$)}

This simple case serves to recap the difference equation method \cite[\S7.5]{linzbook}.
With one spring \eqref{eq:IEfree} is 
a scalar Volterra equation with constant kernel,
so we drop the $j$ index and have
\[
-\sigma(t) - \frac{\scst}{2}\int_0^t \sigma(\tau)d\tau = g(t),
\]
of course equivalent to a first-order ODE IVP.
For an order $p=1$ interpolant we use in each time step the density at the start of that
time step:
\beqs
\sigma(\tau) \approx \sigma^\nu, \qquad t_\nu \le \tau < t_{\nu+1}, \qquad \nu = 1, 2, \dots, \Nt - 1.
\label{interp0}
\eeqs
Substituting this into the above
gives, in terms of the scaled time-step $\al:= \frac{\scst\dt}{2}$, the
$\bigO(\dt)$ accurate explicit time step method
\[
\sigma^{n+1} + \al \sum_{\nu\le n} \sigma^\nu = -g^{n+1}.
\]
Subtracting the equation where $n+1$ is replaced by $n$ shows that the density
obeys the difference equation
\[
\sigma^{n+1} + (\al-1) \sigma^n = -\Delta g^{n+1}
\]
driven by $\Delta g^{n+1}:=g^{n+1}-g^n$.
Its characteristic equation is $p(z) := z + \al-1 = 0$,
and the difference equation is stable if all solutions (roots of $p$) have $|z|<1$
\cite[\S7.5]{linzbook}.
The latter would insure that all homogeneous solutions are exponentially decaying.
This root condition holds whenever $\al \in (0,2)$, thus for all $\dt < 4/\scst$.
(This is the same bound as the forward Euler method for the corresponding first-order ODE.)

Instead, for $p=2$, substituting the $\bigO(\dt^2)$ accurate piecewise-linear interpolant
\beq
\sigma(\tau) \approx \frac{t_{\nu+1}-\tau}{\dt}\sigma^\nu +\frac{\tau-t_\nu}{\dt}\sigma^{\nu+1}
, \qquad t_\nu \le \tau < t_{\nu+1},
\label{interp1}
\eeq
gives the implicit method (analogous to an order-2 linear multistep method or LMM),
\[
\sigma^{n+1} + \frac{\al}{2}\sigma^{n+1}
+ \al \sum_{\nu\le n} \sigma^\nu = g^{n+1}.
\]
Again subtracting the $n+1$ and $n$ cases
gives the difference equation
$
(1+\al/2) \sigma^{n+1} - (1-\al/2) \sigma^n = -\Delta g^{n+1},
$
whose characteristic polynomial has the single root $z = (1-\al/2)/(1+\al/2)$,
inside the unit circle for all $\al>0$. The method is thus unconditionally stable
(analogous to the A-stability of the Crank--Nicholson method for the equivalent ODE
\cite[p.~111]{linzbook}).

\subsection{Stability for two springs ($M=2$)}

Let the spring separation be $L:=|x_2-x_1|$.
The coupled Volterra equations \eqref{eq:IEfree} are then, for $t>0$,
\begin{align}
  \sigma_1(t) + \frac{\scst}{2}\int_0^t \sigma_1(\tau)d\tau
  + \frac{\scst}{2}\int_0^{t-L} \sigma_2(\tau)d\tau
  &= -g_1(t),
  \label{s1}
  \\
  \sigma_2(t)
  + \frac{\scst}{2}\int_0^{t-L} \sigma_1(\tau)d\tau
  +\frac{\scst}{2}\int_0^t \sigma_2(\tau)d\tau
  &= -g_2(t),
  \label{s2}
\end{align}
where we note the symmetry that
the second equation is the first with the labels $j=1$ and $2$ swapped.
The following proves that 1st-order time-stepping is stable up to a maximum $\dt$,
independent of how large $L$ is.

\begin{thm}[Explicit first-order time-stepping for two springs at arbitrary separation]  
  $\qquad$ Let $\dt \in (0,2/\beta)$, and let
  $x_1,x_2\in\R$.
  Then the discrete Volterra equations arising by applying the
  piecewise-constant ($p=1$, i.e., 1st-order)
  interpolant \eqref{interp0} to \eqref{s1}--\eqref{s2} and
  integrating exactly are stable, in the following sense:
  if the right-hand side data $g^n_{1,2}$ are zero for all sufficiently large $n$, the
  densities $\si_{1,2}^n$ remain uniformly bounded, and in fact decay exponentially to
  at worst a constant.
\end{thm}         
\begin{proof}
  Recalling $L:=|x_2-x_1|$,
  let $\kappa := L/\dt$, and let $s = \floor \kappa \in \{0,1,\dots\}$.
The piecewise-constant interpolant integrates exactly over the bisected
time-step $(t_{n-s},t_{n+1}-L)$ to give $(s+1-\kappa) \dt \si^{n-s}$.
Thus the discretized equations are
\beq
\sigma_1^{n+1} +\al \sum_{\nu\le n}\sigma_1^\nu + \tal\si_2^{n-s} + \al\sum_{\nu<n-s}\si_2^\nu
= -g_1^{n+1}
\label{dvie2}
\eeq
plus its $1\leftrightarrow 2$ swapped counterpart,
recalling $\al := \beta \dt/2$,
and defining a bisected time-step weight $\tal := (s+1-\kappa)\al$.
Taking the difference of successive time-steps,
any solution to the above satisfies the coupled difference equations
\beq
\sigma_1^{n+1} + (\al-1) \sigma_1^n + \tal\si_2^{n-s} + (\al-\tal)\si_2^{n-s-1}  = -\Delta g_1^{n+1}
\label{diffeq2}
\eeq
plus its $1\leftrightarrow 2$ swapped counterpart.
Its characteristic equation is a $2\times2$ matrix polynomial,
thus we seek stability by considering the solutions $z\in\C$ of
\beq
\det \left[
\begin{array}{ll}
  z^{s+2} + (\al-1)z^{s+1} &  \tal z+ \al-\tal
  \\
 \tal z+ \al-\tal &     z^{s+2} + (\al-1)z^{s+1}
\end{array}\right]
= 0.
\label{det}
\eeq
This clearly vanishes whenever $z$ is a root of either $p_+$ or $p_-$,
the scalar polynomials
\[
p_{\pm}(z):= z^{s+2} +(\al-1) z^{s+1} \pm (\tal z+ \al-\tal).
\]
Note that $p_-(1) = 0$, giving the factorization
\[
p_-(z) = (z-1)(z^{s+1} + \al z^s 
+\dots+\al z + \al - \tal)
\]
which verifies that $1$ is a simple root of $p_-$, because the remainder
polynomial evaluates to $1+(s+1)\al - \tal > 1+s\al$ when $z=1$.
Also, $p_+(1) = 2\al \neq 0$.
Thus $z=1$ is a simple root of the characteristic equation \eqref{det}.

It remains to prove that all other roots of $p_\pm$, and hence of \eqref{det}, lie strictly
within $|z|<1$.
The following applies to $p_-$ or $p_+$.
We apply Rouch\'e's theorem on the curve $\Gamma_\eps$
defined as the boundary
of the union of the unit disk $|z|<1$ and the overlapping disk $|z-1+\al| < r$, where
$r = \al+\eps$, for $\eps>0$.
We choose a comparison polynomial $q(z) := z^{s+2} + (\al-1)z^{s+1}$,
whose roots always lie in $\Gamma_\eps$
(specifically they are zero with multiplicity $s+1$, and the simple root $1-\al$).
Since $q$ matches the two leading terms of $p_\pm$, only the degree-1 terms are left, so
\[
| p_\pm(z) - q(z) | = |\tal z + \al-\tal| \le \tal|z| + \al-\tal, \quad z\in\C,
\]
since $0<\tal\le\al$. Thus $| p_\pm(z) - q(z) |\le \al(1+\eps)$ for all $z\in\Gamma_\eps$,
since $|z|\le 1+\eps$.
Yet,
\[
|q(z)| = |z^{s+1}| \cdot |z+\al-1|,
\]
and the first term is lower-bounded by $1$ on $\Gamma_\eps$, while the second term
is lower-bounded by $r=\al+\eps$, due to the geometry of the two overlapping disks.
Thus a lower bound for $|q|$ on $\Gamma_\eps$ is $\al + \eps$,
which, by the hypothesis $\al<1$, is strictly larger than the above upper bound on $|p_\pm-q|$.
Thus Rouch\'e implies that $p_\pm$ has $s+2$ roots in $\Gamma_\eps$.
The only point on $|z|=1$ that stays inside $\Gamma_\eps$
in the limit $\eps\to0$ is $z=1$.
Thus the $s+2$ roots must lie in $\{|z|<1\} \cup \{z=1\}$. However, we have already
shown that $z=1$ is either a simple root (case $p_-$) or not a root (case $p_+$).
Thus all other roots lie in $|z|<1$.

Recall that $p_-$ and $p_+$ are the factors of the determinant \eqref{det}.
Then, combining the above, its only root on the unit circle is the simple
root $z=1$, while all $2s+3$ others lie strictly inside the unit circle.
Thus all homogeneous solutions to \eqref{diffeq2} are constant plus exponentially-decaying
components, and thus uniformly bounded.
The same thus holds for solutions to \eqref{dvie2} for all $n$ beyond which the
right-hand sides are zero. \qed
\end{proof}

A few observations on the above somewhat technical proof are warranted:
\begin{enumerate}
\item The difference equation has a larger solution space
than the discrete integral equation: it allows a neutrally stable mode
(but we do not know whether it is excited in the scattering setting).
\item Attempting to apply Rouch\'e's theorem as usual on a circle $|z|=r$, taking
$r\to1^+$, fails to exclude the possibility of other roots on $|z|=1$ which
may have arbitrarily high multiplicity.
This explains the need for the novel contour used.
\item Applying the symmetrization argument used below
  unfortunately forces a pessimistic stability bound of the form
  $\dt < \bigO(1/\kappa)$, much worse than achieved above via polynomial roots.
\end{enumerate}

We next prove that 2nd-order time-stepping is stable no matter how large $\dt$ is, for two
``nearby'' springs in the sense of Section~\ref{sec:local},
meaning $L<\dt$ (so the full implicit nature of the scheme is engaged, and the matrix $S$ in
\eqref{eq:linearSystem} is nonzero).
For this to hold there is a mild upper bound on $L$.

\begin{thm}[Implicit second-order time-stepping for two springs separated by less than $\dt$]  
  \label{thm:stability1}
Let $L <2/\beta$, and let $\dt>L$.
Let $\gvec := [g_1^{1}, \dots, g_1^{\Nt}, g_2^{1}, \dots, g_2^{\Nt}]$ represent the data vector
for $\Nt$ time steps.
Let $\sigvec := [\sigma_1^{1}, \dots, \sigma_1^{\Nt}, \sigma_2^{1}, \dots, \sigma_2^{\Nt}]$ be the
vector of densities computed using the discrete Volterra equations arising by applying the linear ($p=2$) interpolant \eqref{interp1} to \eqref{s1}--\eqref{s2}, and integrating exactly.
Then there is a constant $C>1$, independent of the data, such that
\beq
\norm{\sigvec}\leq C \norm\gvec,
\label{sigbnd}
\eeq
where $\norm{\cdot}$ is the 2-norm on $\mathbb R^{2\Nt}$, i.e., the scheme is stable.
\end{thm} 

\begin{proof}
Using $\kappa := L/\dt$, $\kappa <1$, the linear interpolant~\eqref{interp1} integrates exactly over the bisected time-step $(t_n, t_{n+1} - L)$ to give $\frac{(1-\kappa)}{2}[(1 + \kappa)\sigma^n + (1 - \kappa) \sigma^{n+1}]$.
The discretized integral equations become 
\beq
\left(1 + \frac{\al}{2}\right)\sigma_1^{n+1}+\al \sum_{\nu\le n}\sigma_1^\nu + \al \sum_{\nu<n}\sigma_2^\nu + \xi_1 \sigma_2^n  + \xi_2 \sigma_2^{n+1} 
= -g_1^{n+1},
\label{dvie3}
\eeq
plus its $1\leftrightarrow 2$ swapped counterpart, defining $\xi_1: = \alpha(2 - \kappa^2)/2$, and $\xi_2: = \alpha(1 - \kappa)^2/2$, again with $\al=\beta \dt/2$.
The scheme in \eqref{dvie3} can be expressed as
\beq\label{eq:stability_system}
\left( a I+A\right) \sigvec = - \gvec,
\eeq
where $a := 1 + \al/2$, and $A$ is the $2\Nt\times2\Nt$ matrix
with lower-triangular blocks given by
\beqs
A = \begin{bmatrix}
L_0&\tilde L_0\\
\tilde L_0 & L_0
\end{bmatrix}, \quad L_0 = \begin{bmatrix}
0&\dots&\dots& 0\\
\al&\ddots&& \vdots\\
\vdots &\ddots &\ddots&\vdots\\
\al&\dots  &\al &0
\end{bmatrix}, \quad
\tilde L_0 = 
\begin{bmatrix}
\xi_2 & 0 & \dots &\dots& 0\\
\xi_1 & \ddots &\ddots &&\vdots\\
\al & \ddots & \ddots &\ddots &\vdots\\
\vdots & \ddots & \ddots & \ddots &0 \\
\al & \dots &\al & \xi_1 & \xi_2
\end{bmatrix}.
\eeqs

We adapt a symmetrization argument from \cite[Sec.~4]{explicitheat19}
to show that $\norm{\sigvec}$ is bounded by the data $\norm{\gvec}$.
We multiply Equation~\eqref{eq:stability_system} by $\sigvec^{T}$ from the left
\beqs\label{eq:stability_normM1}
a \norm{\sigvec}^2 + \sigvec^{T}A\sigvec = -\sigvec^{T}\gvec, 
\eeqs
and note that for $W = A + A^{T}$, $\sigvec^{T}A\sigvec = \frac{1}{2}\sigvec^{T}W\sigvec$,
so that
\beq\label{eq:stability_normW}
a\norm{\sigvec}^2 + \frac{1}{2}\sigvec^{T}W\sigvec = -\sigvec^{T}\gvec. 
\eeq
We further partition $W$ into a sum of two matrices $V$ and $Q$ by introducing a positive parameter $\theta$ and writing $W = V + Q$. The matrix $Q$ takes the form
\beq
{Q = \al \bm{1}\bm{1}^{T} - (\al + \theta) I},\eeq
where $\bm{1}$ is the $2\Nt\times 1$ vector of ones, and $I$ is the $2\Nt\times 2\Nt$ identity matrix. The matrix $V$ is
\beq
V = \begin{bmatrix}
\theta I & V_1\\
V_1 & \theta I
\end{bmatrix}, 
\eeq
where $V_1$ is
tridiagonal with $\omega_2 := 2\xi_2 - \al = \al\kappa(\kappa - 2)$ on the main diagonal and ${\omega_1 := \xi_1 - \al = -\kappa^2\al/2}$ on the sub-diagonals.
The spectral norm $\|V_1\|$ at least $|\omega_2+2\omega_1| = 2\kappa\al$ because this is
the eigenvalue for the constant vector, and by Gershgorin's theorem this is also
the upper bound, so $\|V_1\|=2\kappa\al$.
The block structure of $V$ means that it is positive semi-definite if
\beq
\theta \ge \|V_1\| = 2\kappa\al = L\beta,
\label{theta}
\eeq
using the definitions of $\kappa$ and $\al$.
With the partitioning $W=V+Q$, equation \eqref{eq:stability_normW} yields
\beqs
a\norm{\sigvec}^2 + \frac{1}{2}\sigvec^{T}\left(V + \al\ones - \al I - \theta I\right)\sigvec = -\sigvec^{T}\gvec,
\eeqs
but we may drop the nonegative quantities $\sigvec^{T}\ones\sigvec\geq 0$ and $\sigvec^{T}V\sigvec\geq 0$, and bound $-\sigvec^{T}\gvec\leq \norm{\sigvec}\norm{\gvec}$ by Cauchy--Schwarz,
to get
\beq\label{eq:stability_sigmaBound}
\left(1 - \frac{\theta}{2}\right)\norm{\sigvec} \leq \norm \gvec.
\eeq
This bound is only useful if $\theta<2$; by the hypothesis $L\beta<2$ this is
compatible with \eqref{theta}, and one may pick $\theta = L\beta$.
The constant in \eqref{sigbnd} may thus be chosen as $C = 1/(1-L\beta/2)$.
\qed
\end{proof}

\subsection{A convergence result for two springs}

Our final result is convergence of the 2nd-order accurate time-stepping scheme for two springs,
in the style of \cite[Thm.~7.2]{linzbook}
(this theorem cannot be applied directly because of the wave propagation delays in our kernels).
Given a time step $\dt$, let $\sigma_1^n$, $\sigma_2^n$
be the solutions computed using the discrete Volterra equations arising from substituting the linear ($p=2$)
interpolant \eqref{interp1} in \eqref{s1}--\eqref{s2}, and integrating exactly. Let $\sigma_1(t_n)$, $\sigma_2(t_n)$ be samples of the true solutions of the Volterra equations in
\eqref{eq:IEfree}. Define the error vector at time $t_n=n\dt$ to be $\errvec{n} := [\err{1}{n},\err{2}{n}]^{T}$ where the error at each scatterer is $\err{j}{n} := \sigma_j(t_n) - \sigma_j^n$, $j = 1, 2$.
 
\begin{thm}\label{thm:convergence}
Let $0<\dt < \min(2/\beta, L)$, and consider the
discrete solution error vector defined above, over a fixed solution interval $t\in(0,T)$.
Suppose that the densities initially vanish, i.e., $\sigma_{1,2}(0) = 0$.
Then, using $\norm{\cdot}_1$ as the 1-norm on $\mathbb R^{2}$,
\beq
\norm{\errvec{n}}_1 = \Oh(\dt^2) \quad \text{as} \quad \dt\rightarrow 0,
\quad \mbox{ uniformly over } \; 0 \leq n \leq \Nt = \lceil T/\dt \rceil.
\eeq
\end{thm}

\begin{proof}
Again set $\kappa = L/\dt$, and $s = \floor\kappa\geq 1$.
Integrating the linear interpolant in~\eqref{interp1} over the bisected time-step $(t_{n - s}, t_{n+1} - L)$ gives 
\beq
\frac{\dt}{2}(1 - \kappa + s)(1 + \kappa - s)\sigma^{n-s} + \frac{\dt}{2} (1 - \kappa + s)^2\sigma^{n - s + 1}, 
\eeq 
so that the discrete equations become
\beq\label{dvie4}
\left(1 + \frac{\alpha}{2}\right) \sigma_1^{n+1} + \alpha\sum_{\nu\leq n}\sigma_1^{\nu}  + \alpha\sum_{\nu<n - s} \sigma_2^{\nu} + \xi_1\sigma_2^{n - s} + \xi_2\sigma_2^{n+1 - s} = - g_1^{n+1},
\eeq 
plus its $1\leftrightarrow2$ swapped counterpart, where 
\beq
\begin{split}
\xi_1 &:= \frac{\alpha}{2}\left[(1 - \kappa+ s)(1 + \kappa- s) + 1\right], \\
\xi_2 &:= \frac{\alpha}{2}\left(1 - \kappa + s\right)^2. 
\end{split}
\eeq
Subtracting \eqref{dvie4} from the true integral equation \eqref{eq:IEfree} evaluated at $t_{n+1}$ gives 
\beq\label{eq:errorEq}
\begin{split}
\err{1}{n+1} = &-\frac{\beta}{2}\left(\int_{0}^{t_{n+1}}\sigma_1(\tau)d\tau - \frac{\dt}{2}\sigma_1^{n+1} - \dt\sum_{\nu \leq n}\sigma_1^{\nu}\right)\\
 &- \frac{\beta}{2}\left(\int_{0}^{t_{n+1}-L}\sigma_2(\tau)d\tau -\dt\sum_{\nu < n-s}\sigma_2^{\nu} - \tilde\xi_1\sigma_2^{n - s} - \tilde\xi_2\sigma_2^{n - s + 1}\right),
\end{split}
\eeq
and its $1\leftrightarrow2$ swapped counterpart, where $\tilde\xi_\ell: = (2/\beta)\xi_\ell$, $\ell = 1, 2$. 
Define the local consistency errors
\beq
\begin{split}
\con{1}{n+1} &= \int_{0}^{t_{n+1}}\sigma_1(\tau)d\tau - \frac{\dt}{2}\sigma_1(t_{n+1}) - \dt\sum_{\nu \leq n}\sigma_1(t_\nu), \\
\tcon{2}{n+1} &= \int_{0}^{t_{n+1} - L}\sigma_2(\tau)d\tau - \dt\sum_{\nu <n - s}\sigma_2(t_{\nu}) - \tilde\xi_1\sigma_2(t_{n - s}) - \tilde\xi_2 \sigma_2(t_{n - s + 1}), \\
\end{split}
\eeq
along with the $1\leftrightarrow 2$ swapped counterparts. 
By design of the 2nd-order accurate quadrature scheme,
noting that no endpoint correction is needed at $t=0$ due to the densities
vanishing there, we have
\beq
\max_{0\leq n\leq \Nt}\abs{\con{1}{n}}\leq c\dt^2, \qquad \max_{0\leq n\leq \Nt}\abs{\tcon{2}{n}}\leq \tilde c\dt^2,
\eeq
where $c$ and $\tilde c$ are constants. 
The error equation in \eqref{eq:errorEq} becomes
\beq
\begin{split}
\err{1}{n+1} = & -\frac{\beta}{2}\left(\con{1}{n+1} + \frac{\dt}{2}\err{1}{n+1} + \dt\sum_{\nu \leq n}\err{1}{\nu}\right)\\
& -\frac{\beta}{2}\left(\tcon{2}{n+1} + \dt\sum_{\nu <n-s}\err{2}{\nu} + \tilde\xi_1\err{2}{n - s} +\tilde\xi_2\err{2}{n - s + 1}\right).
\end{split}
\eeq
Noting that $|\tilde\xi_1|, |\tilde\xi_2| \le \dt$, we write 
\beq
\left(1 - \frac{\beta\dt}{4}\right)\abs{\err{1}{n+1}}\leq \frac{\beta}{2} \left(\abs{\con{1}{n+1}} + \abs{\tcon{2}{n+1}}\right) + \frac{\dt\beta}{2}\sum_{\nu \leq n}\left(\abs{\err{1}{\nu}} + \abs{\err{2}{\nu}}\right).
\eeq
We then have
\beq\label{eq:errBound1}
\left(1 - \frac{\beta\dt}{4}\right)\norm{\errvec{n+1}}_1 \leq \frac{\beta}{2} \left(\abs{\con{1}{n+1}} + \abs{\con{2}{n+1}} + \abs{\tcon{1}{n+1}} + \abs{\tcon{2}{n+1}}\right) + \dt\beta\sum_{\nu \leq n} \norm{\errvec{\nu}}_1.
\eeq
By consistency, there exists a constant $C$ such that 
\beqs
\abs{\con{1}{n+1}} + \abs{\con{2}{n+1}} + \abs{\tcon{1}{n+1}} + \abs{\tcon{2}{n+1}}\leq C\dt^{2}, \quad \forall n = 0, \dots, \Nt-1.
\eeqs
Provided that $\beta\dt<4$, we may write the inequality in~\eqref{eq:errBound1} as
\beq
\norm{\errvec{n+1}}_1\leq \frac{\beta}{2\left(1 - \frac{\beta\dt}{4}\right)}C\dt^2 + \frac{\dt\beta}{\left(1 - \frac{\beta\dt}{4}\right)}\sum_{\nu \leq n}\norm{\errvec{\nu}}_1. 
\eeq
Applying the discrete Gronwall inequality \cite[Thm.~7.1]{linzbook}, we have 
\beq
\norm{\errvec{n}}_1\leq \frac{\beta}{2\left(1 - \frac{\beta\dt}{4}\right)}C\dt^2\exp\left(\frac{\beta \dt n}{1 - \frac{\beta\dt}{4}}\right). 
\eeq
Specifically, for $\beta\dt<2$, since $n\dt<T$, we arrive at the following estimate
\beq
\norm{\errvec{n}}_1 \leq \beta C \dt^2 e^{2\beta T}.
\eeq
This proves the claim.
\qed
\end{proof}

\section{Numerical tests and examples}
\label{sec:results}
We demonstrate the high-order convergence and computational performance of the WFP method using numerical examples. In section~\ref{sec:manufacturedSolution}, we use the method of manufactured solutions to generate convergence results, while section~\ref{sec:timing} investigates the computational performance of the scheme for free-space spring scattering problems.
We investigate physical properties of scattering from randomly placed springs in
Section~\ref{sec:randomSprings}, demonstrating wave localization.
We finish with demonstrations of narrow-band (high-Q) filtering by systems of
equi-spaced springs in Section~\ref{sec:uniformSprings}.

In all experiments we choose window tolerance parameter $\epsilon = 10^{-12}$, and Nyquist fraction $\gamma = 0.5$, as described in section~\ref{sec:periodic}.
For nonuniform FFTs we use the FINUFFT library \cite{finufft,finufftlib},
which exploits shared-memory parallelism.
Our implementation is in MATLAB (using version R2023b),
and the local evaluation is done by precomputing a sparse matrix
which acts on the stack of the recent $\mmax$ density vectors.
Tests are run on an Apple M2 Max laptop (released 2023) with 64 GB RAM,
apart from the scaling test of section~\ref{sec:timing} which is
run on a Ice Lake node with two 32-core Xeon Platinum 8362 2.8 GHz CPUs
(released 2021) and 1024 GB RAM.

\subsection{Convergence test}\label{sec:manufacturedSolution}

To assess the convergence of the WFP method
described in sections \ref{sec:periodic}--\ref{s:method}
we use the method of manufactured solutions:
we generate a free-space scattered potential $u_{\text{ex}}$
from analytic densities,
read off its data $g_j(t)$, apply the WFP method to solve
numerically the free-space Volterra integral equations \eqref{eq:IEfree}
with this data,
then from the resulting solution densities evaluate a
new scattered potential $u$.
We report the error $u - u_{\text{ex}}$.
The point is to compare to an analytic solution,
noting that no analytic solution is known for actual scattering problems.

We place the scatterers uniformly at random at locations $x_j\in[-1,1]$,
enforcing a minimum separation of $10^{-4}$ between any two scatterers.
As densities we choose Gaussian pulses
$$\sigma_j(t) = e^{-\mu_j(t - t_{0,j})^2}, \qquad j = 1, \dots, M,$$
with random values assigned to the inverse variances $\mu_j\in[40,50]$
and peak times $t_{0,j}\in[1,3]$.
These densities generate
$u_{\text{ex}} = {\cal S} [\Gamma,\vec{\sigma}]$
according to the representation \eqref{eq:singleLayerSol}.
We assign random spring strengths $\scst_j\in[0.1,3]$, for $j = 1,\dots, M$.
The data $g_j(t), j = 1,\dots, M,$ is then extracted as
\beq
g_j(t) = -\sigma_j(t) - \scst u_{\text{ex}}(x_j,t),
\eeq
recalling \eqref{kink} and the jump relation \eqref{JR}.

\begin{figure}[th]
	\centering
	\includegraphics[width=7.5cm]{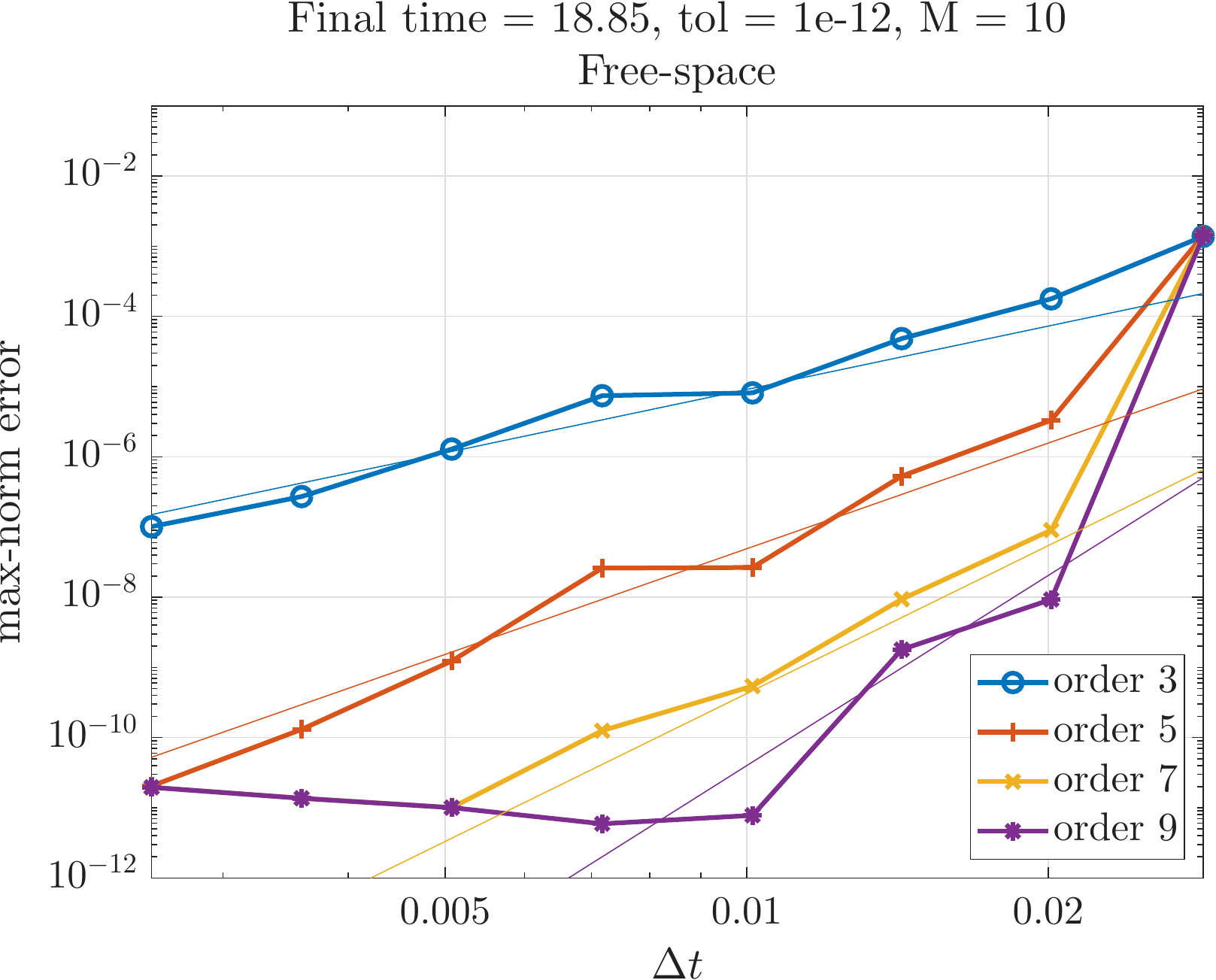}
	\includegraphics[width=7.5cm]{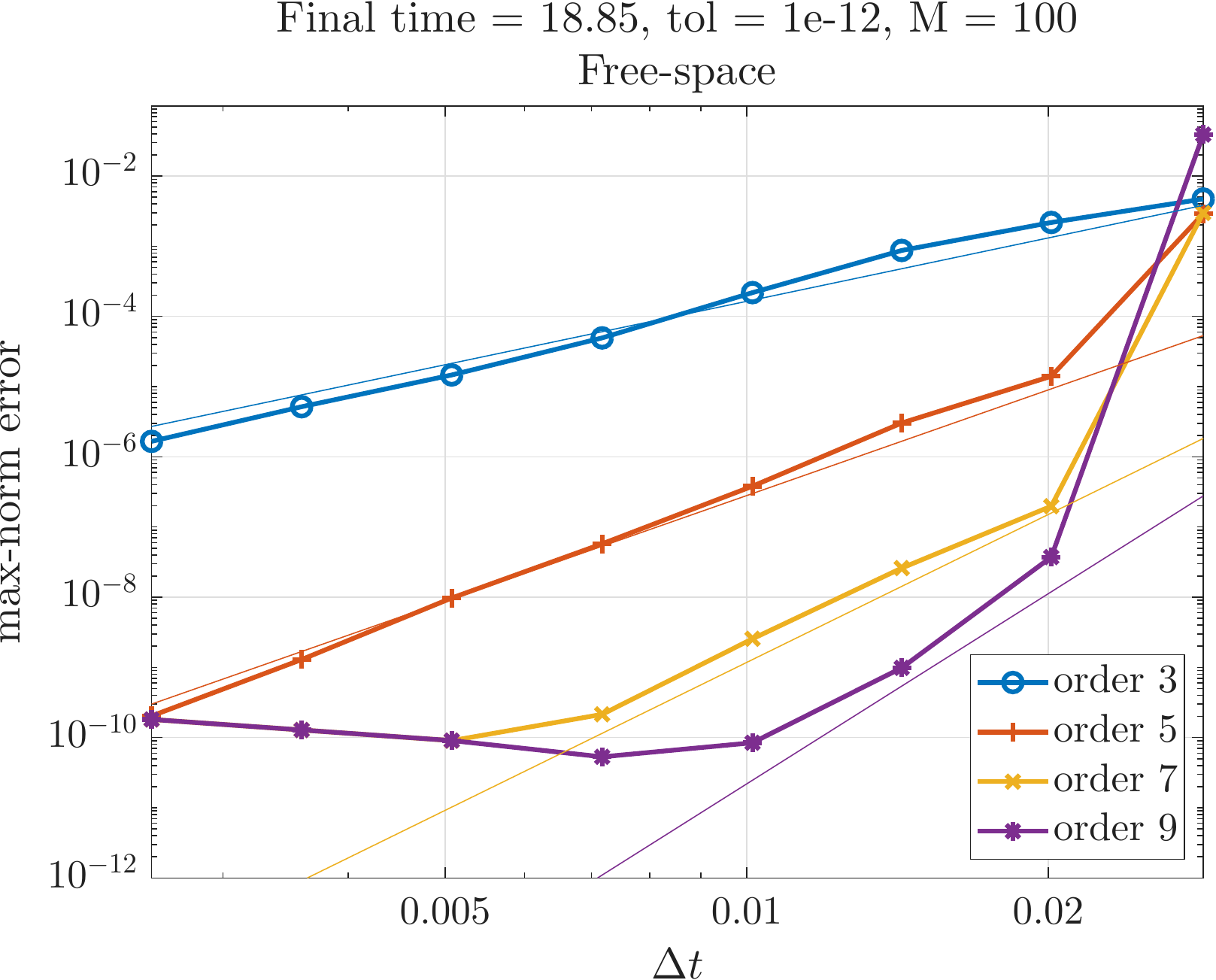}
	\caption{Convergence to a manufactured solution, for the WFP method at interpolation orders $p=2,4,6,8$. The observed orders $p+1$ are shown by the thin reference lines.
Left plot shows $M=10$ springs; right plot shows $M=100$ springs.
Convergence bottoms out near the window tolerance $\eps=10^{-12}$.
See section~\ref{sec:manufacturedSolution} for details.}
        \label{fig:msConvergence}
\end{figure}

In Figure~\ref{fig:msConvergence} we use this to test the WFP method convergence
with respect to $\dt$, at $p=2, 4, 6, 8$ interpolation nodes
as defined in Section~\ref{sec:local}.
The maximum error between scattered field
$u$ and $u_{\text{ex}}$ is
measured over a uniform space-time grid of size $10\times10$,
covering $0<t<T=6\pi$.
We observe errors $\bigO((\dt)^{p+1})$, that is,
one order higher than the order-$p$ polynomial
interpolation in Section~\ref{sec:local}.
The additional $\dt$ factor appears in the local quadrature error 
since we are making an $\bigO((\dt)^{p})$ error on an interval of size $\bigO(\dt)$.
In our Volterra integral equation, 
this local quadrature error dominates the total error since the history part is computed with spectral accuracy.
As expected, this algebraic
convergence stops at the window tolerance $\bigO(\eps)$,
and we note that the scaling with the number of scatterers is well
explained by $M\eps$.

\subsection{Scaling of CPU time with problem size}\label{sec:timing}

Here we measure the computational performance of WFP
for large numbers of scatterers $M$.
We choose sources randomly located in $[-2,2]$ with minimum separation of $10^{-6}$ between any two sources, and random spring constants $\beta_j\in[0.1,3]$.
In this and the following sections we
solve the free-space scattering problem with Gaussian incident pulse
\beq\label{eq:incidentPulse}
f(t) = e^{-\mu(t-t_0)^2},
\eeq
giving the incident wave \eqref{uinc}.
In this section we set $\mu=30$ and $t_0=-3$.
For each $M$, we select $\dt$ such that $\Mtyp \approx 100$ (see section \ref{s:compfull}).
We set $p=6$ interpolation nodes (7th-order scheme).
This typically gives an accuracy of 10 digits at smaller $M$, dropping
to 6 digits at $M=10^{6}$.
We report the CPU time required per time-step,
and in this section run experiments on the Xeon node described above.
This does not include (history and local) evaluation at target points, which becomes negligible at large $M$. 
We do not report precomputation time, since
it is less than $100$ times the cost of one time-step,
and in all of our experiments this is negligible compared to the total cost.

\begin{figure}[th]
	\centering
	\includegraphics[width=8cm]{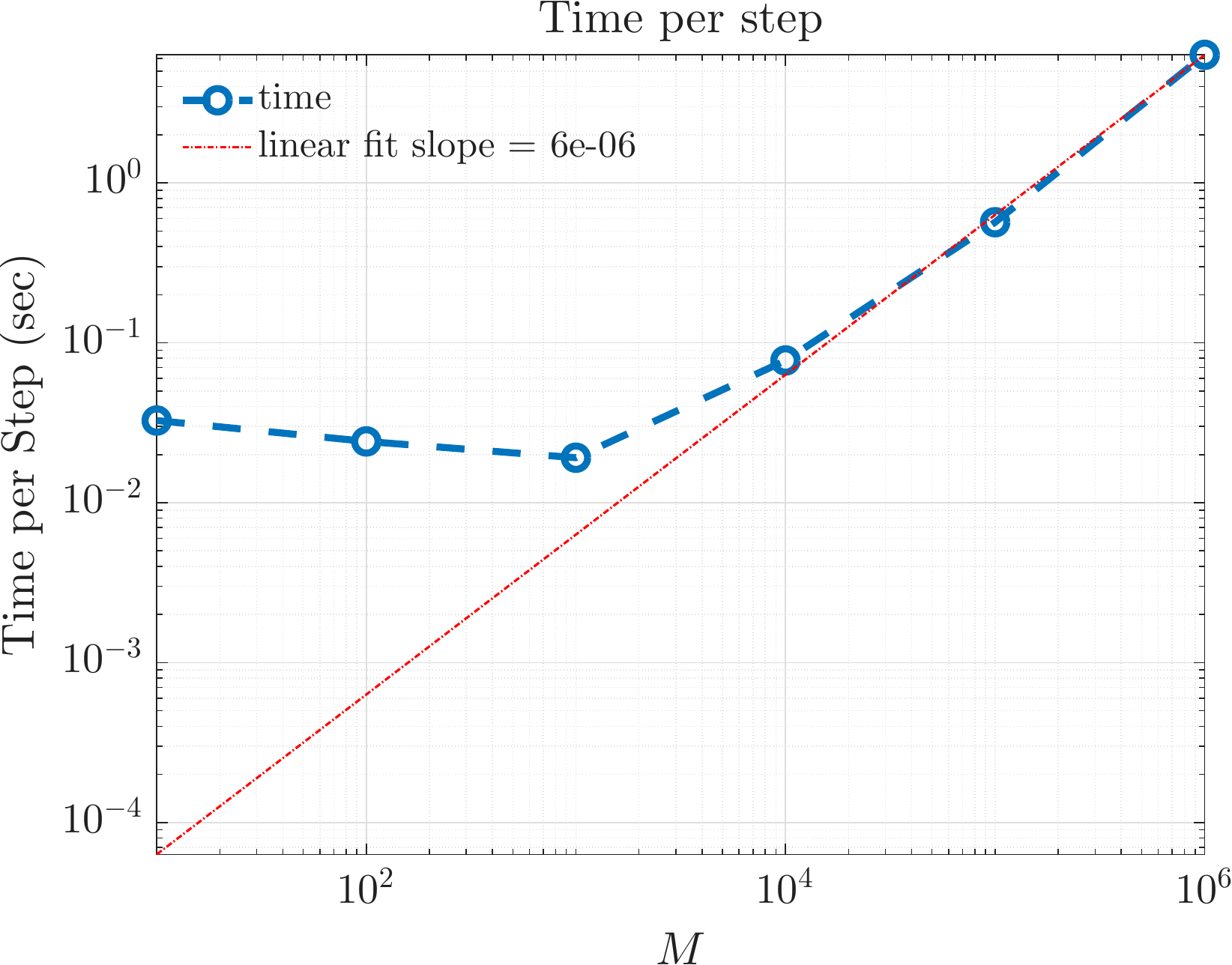}
	\caption{Computational cost per time-step, on a Xeon compute node, for the 7th-order accurate implementation of the WFP scheme as a function of $M$, the number of spring scatterers. See section~\ref{sec:timing} for details.}
	\label{fig:timing}
\end{figure}

Figure~\ref{fig:timing} shows the results: they are
consistent with a cost asymptotically linear in $M$,
approaching a throughput of about 150000 scatterer-time-steps per second.
The time-step is dominated by a sparse matrix-vector multiplication
and a pair of NUFFTs.
Note that since $K$ for the Fourier series grows like $1/\dt$, which
grows like $M$ in order to keep $\Mtyp$ fixed, the NUFFTs should
cause the complexity to be dominated by $\bigO(M \log M)$ work per time-step.
However, it appears that even at the largest $M=10^6$ tested, this
dominance is not yet visible.
The storage is relatively mild: the example with a million scatterers
needs only around 2 GB RAM.

\begin{figure}[th]
  \centering
        \foreach\scstVal in {3,5,10}{
	\includegraphics[width=4cm]{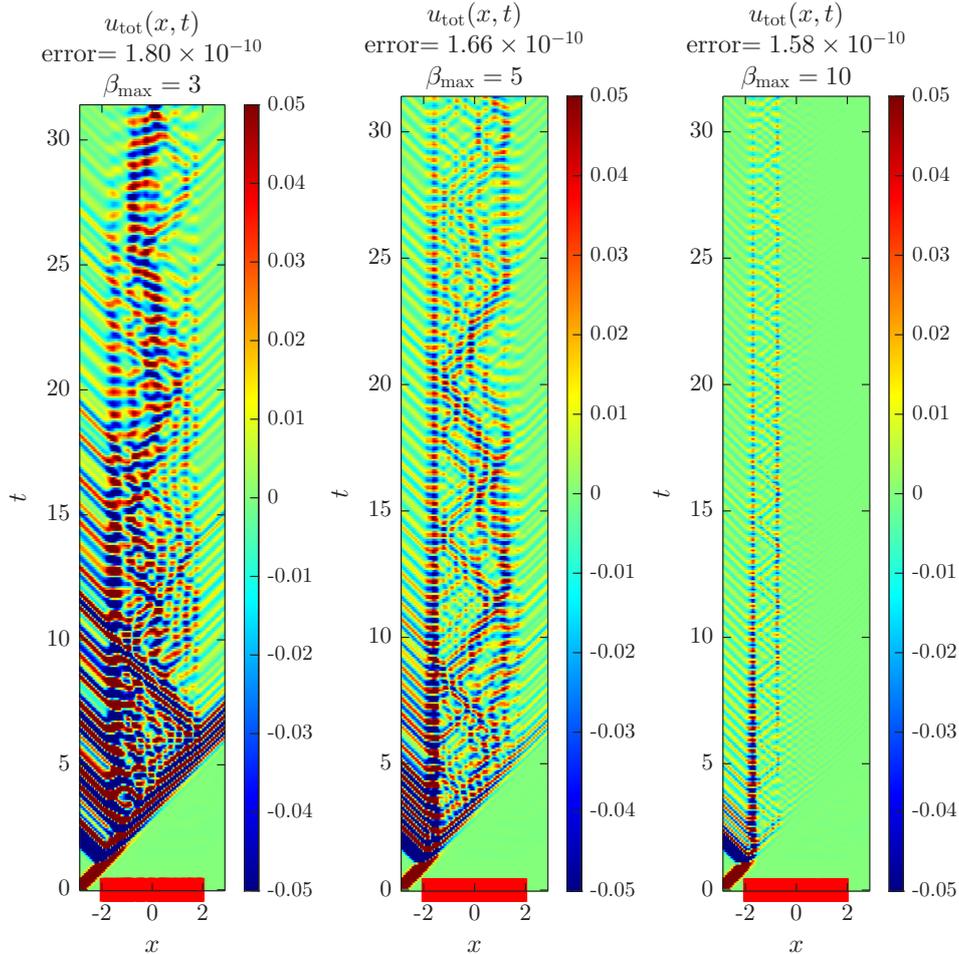}
	}
	\caption{9th-order accurate ($p=8$) computed total field when scattering from a region containing $M = 150$ randomly located springs (appearing as a dense red block), with random strengths $\scst_j\in[0.1,3]$ (left), $\scst_j\in[0.1,5]$ (middle), or $\scst_j\in[0.1,10]$ (right). The time step is $\dt = 3.703\times10^{-3}$.
This takes about 40 seconds of CPU time to compute
(of which 0.6 is precomputation), using $\Nt = 8484$ time-steps.}
	\label{fig:localizationM150}
\end{figure}

\subsection{Scattering in 1D disordered media}\label{sec:randomSprings}

It is a famous result that a spatially-varying random potential
can induce localization of quantum waves \cite{anderson58}.
Similar effects occur
in acoustic and photonic (wave equation) systems in random media.
The 1D case, where localization is strongest, is of
interest in physics \cite{vib1983,condat86,Disordered2009}.
In 1D all frequencies 
are expected to localize, meaning the absence of
propagation in an unbounded random medium, and trapping at sites.
The amplitude decays exponentially away from such sites, on average,
defining a ``localization length'', which may depend
on frequency \cite{Freilikher92}.

In this section, we study scattering from large numbers of
springs at random locations to
explore wave localization, and to demonstrate our solver's wide-band capability
in challenging problems.
We use the WFP method to
solve numerically the free-space scattering problem from the introduction,
with the Gaussian incident pulse \eqref{uinc} defined in~\eqref{eq:incidentPulse}. We take $x_j\in[-2,2]$ to be randomly spaced, such that the minimum distance between any two scatterers is $10^{-4}$.
For additional disorder we choose the spring constants $\beta_j$ randomly.

Figure~\ref{fig:localizationM150} plots the computed total field $\utot = u + \uin$ with space horizontally and time vertically (a ``space-time'' diagram),
computed on a $400\times400$ regular grid,
for $M=150$ scatterers.
Three different random ranges of spring strength are tested (see caption).
The pulse parameters in \eqref{eq:incidentPulse} were $\mu = 30$
and $t_0 = -3$.
We fix the time step $\dt$ such that $\Mtyp \approx 10$,
and run the simulation up to a final time of $T=10\pi$.
The compute time (see caption) is very modest, less than a minute.
A self-convergence study is performed
(comparing to a simulation using $\dt/2$),
giving the estimated uniform errors stated at the
top of each plot: the solutions have around 10 digits accuracy.

We observe that the high frequency part of the incident pulse
passes through with little delay in the left and middle cases in
Figure~\ref{fig:localizationM150}, whereas lower frequency
waves are reflected, or trapped
for long times inside the spring region.
Clear resonant trapping at sites is visible in the right-most case.
Here, almost all of the incident wave is reflected at the first bounce.
This can be explained as follows, without invoking localization.

\begin{remark}[Homogenization theory for frequency filtering]
\label{r:KleinG}
We apply a simple homogenization argument to the wave equation
\eqref{WE} with a dense random set of springs \eqref{kink},
whose mean separation is $\Delta x$
and mean strength $\overline{\scst}$.
Replacing the springs with a distributed spring
with the same mean restoring force per unit length
gives the Klein--Gordon equation
\beq
\pt^2 u - \px^2 u + \om_0^2 u = 0, \qquad \mbox{ where } \;
\om_0^2 = \overline{\scst} / \Delta x,
\eeq
as an effective PDE for the region $[-2,2]$.
Here $\om_0$ is the cut-off frequency below which waves cannot propagate
(the dispersion relation is $\om = \sqrt{\om_0^2 + \kappa^2}$).
Applying this to the three strength distributions in
Figure~\ref{fig:localizationM150}
gives $\om_0 \approx 7.6$ for the first case, rising to
$\om_0 \approx 13.8$ for the third case.
The incident frequency content is the Gaussian
$\hat{f}(\om) \propto e^{-\om^2/4\mu}$, which has standard deviation
$\sqrt{2\mu} \approx 7.7$.
Thus in the first case a significant portion of the incident
frequencies lie above $\om_0$, explaining the propagation of
high frequencies right through the spring region.
In the last case $\om_0$ is nearly at two standard deviations,
explaining why very little energy propagates through the medium.
\end{remark}

\begin{figure}[th]
	\centering
	\includegraphics[width=7.5cm]{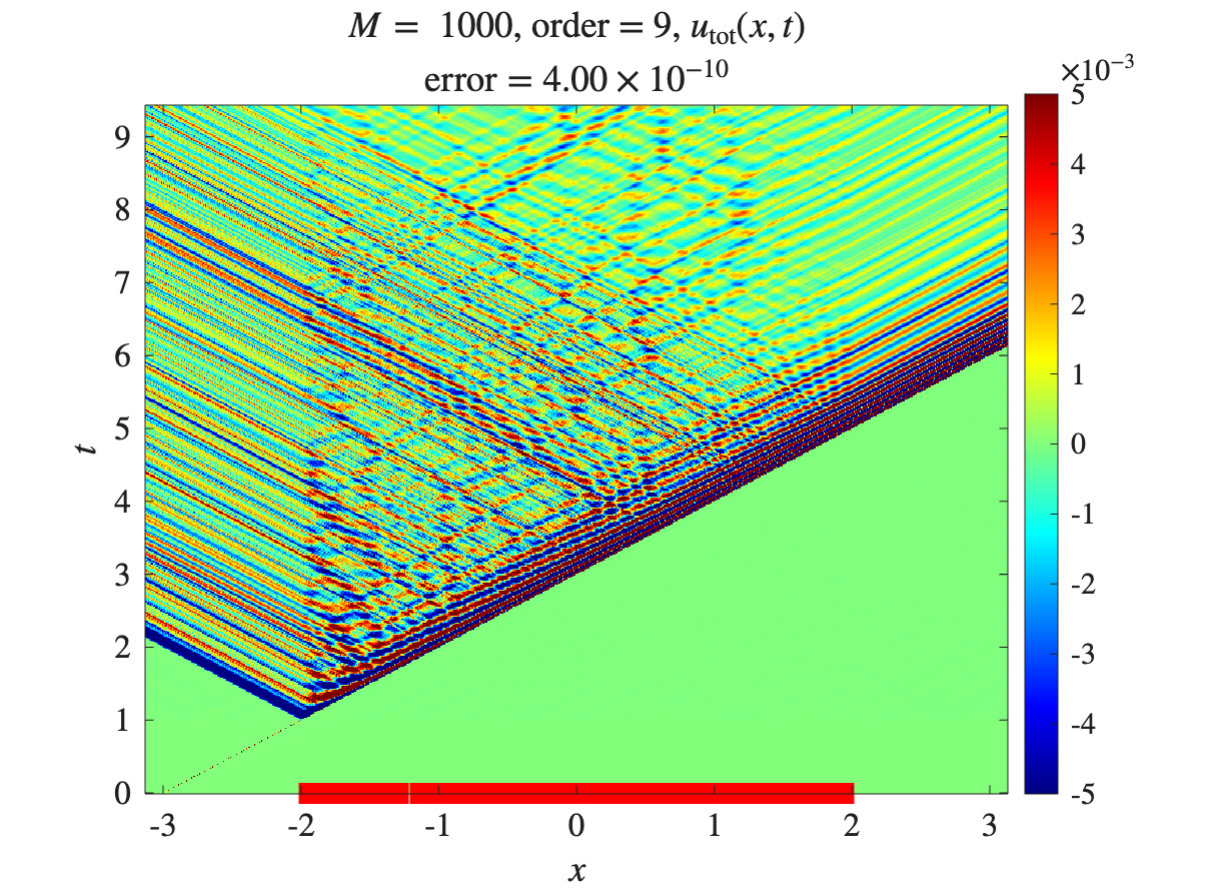}
	\includegraphics[width=7.5cm]{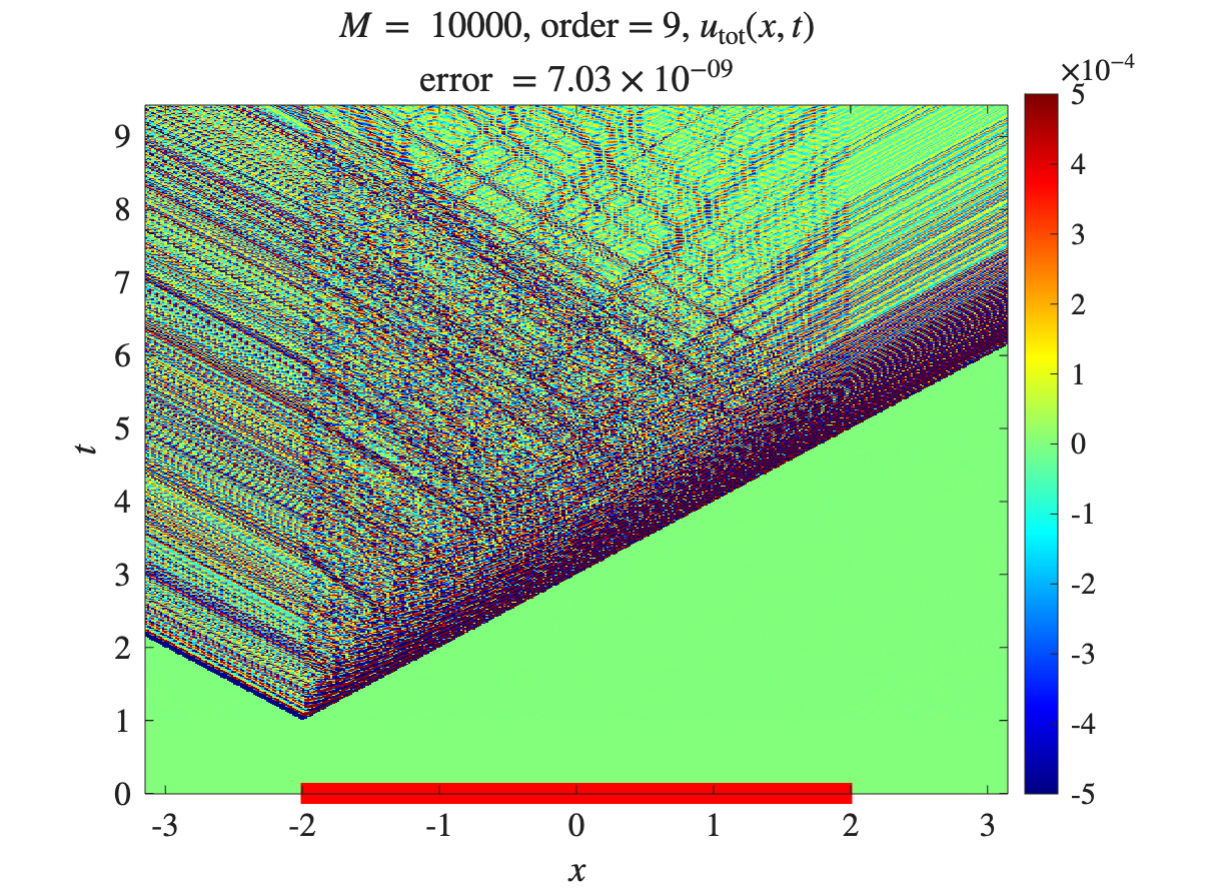}
	\caption{Ninth-order accurate computed total field for a narrow incident wave with $\mu = M^2$ and $M =10^3$ (left) or $M=10^4$ (right) random spring scatterers. The incident wave is almost too narrow to see in either plot.
       For $M = 10^3$, $\dt = 1.308\times10^{-4}$, and the runtime is
       24 minutes (1.2 sec for precomputation), and 0.02 sec per time-step for
       $\Nt = 72044$ total time steps.
       For $M = 10^4$, $\dt = 1.308\times10^{-5}$, and the runtime
       is 29.5 hours (9 sec precomputation),
       at 0.15 sec per time-step for $\Nt = 720440$ time steps.}
	\label{fig:largeMscattering2}  
\end{figure}

Figure~\ref{fig:largeMscattering2}
shows challenging examples with large numbers of random springs: $M=10^3$ (left)
and $M=10^4$ (right).
As $M$ grows, $\dt$ must decrease in order that the number of local neighbors
$\Mtyp$ remain $\bigO(1)$ to keep the method linear-scaling
(in other words, a longer $\dt$ would not be cheaper, due to the local cost),
but also for stability reasons.
However, this means that the scheme can accurately handle a growing
incident wave bandlimit of $\bigO(M)$.
We thus choose the incident pulse to have a width
$\bigO(1/M) = \bigO(\Delta x)$, by setting $\mu = M^2$.
Regardless of $M$,
such a pulse has a standard deviation (width) $\sqrt{8}$ times the
mean inter-spring spacing $\Delta x = 4/M$.
As before, $t_0=-3$.
The spring constants $\scst_j$ are chosen randomly from the interval $[0.1,3]$, and $\dt$ is such that the typical number of neighbors for each source is $\Mtyp = 50$ for $M = 10^3$, and $\Mtyp = 500$ for $M = 10^4$.

We observe complex multiple reflections, dispersion (of a strong transmitted wave), and trapping; a full study is beyond this paper.
Since $\bigO(M)$ time-steps are needed to simulate a couple of wave passage
times, the total simulation cost grows like $\bigO(M^2 \log M)$.
This is apparent since growing $M$ by a factor of 10 turns a 24-minute
run into a 28-hour run.
By comparing against runs at $\dt/2$, the estimated errors
are over nine digts for $M=10^3$, and over eight digits for $M=10^4$.

\begin{figure}[th]
	\centering
	\includegraphics[width=7.5cm]{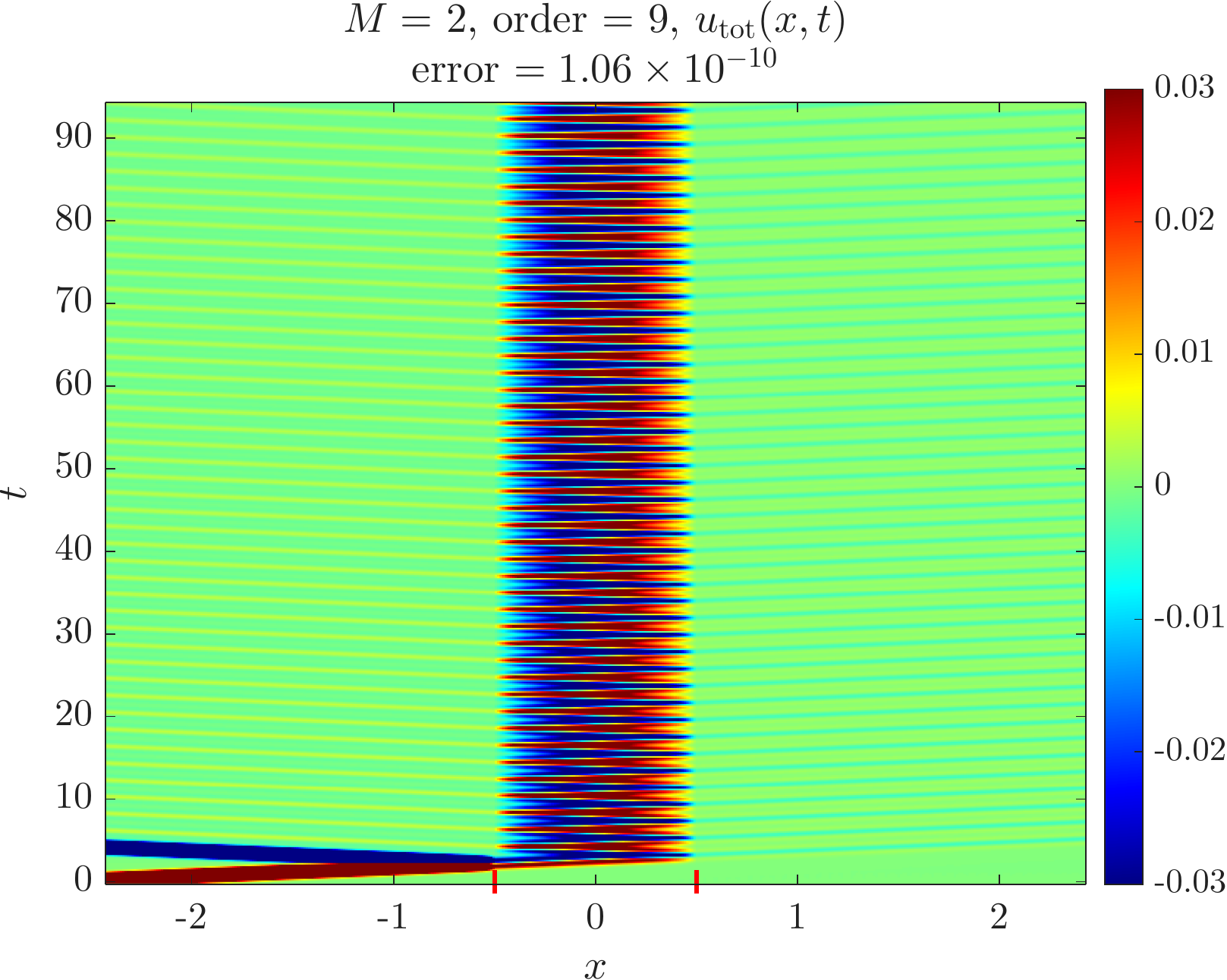}
	\caption{Scattering from two springs (indicated by red tick marks on the $x$-axis) a unit distance apart, with spring constants $\beta = 100$.
Note that most of the incident energy is reflected; what passes through
is very narrow in frequency.
This simple simulation takes 16 seconds, with $\Nt=9426$ time steps of
$\dt=0.01$. See Section~\ref{sec:uniformSprings}.}
	\label{fig:FP_heatmap} 
\end{figure}

\begin{figure}[th]
	\centering
        \includegraphics[width=7.5cm]{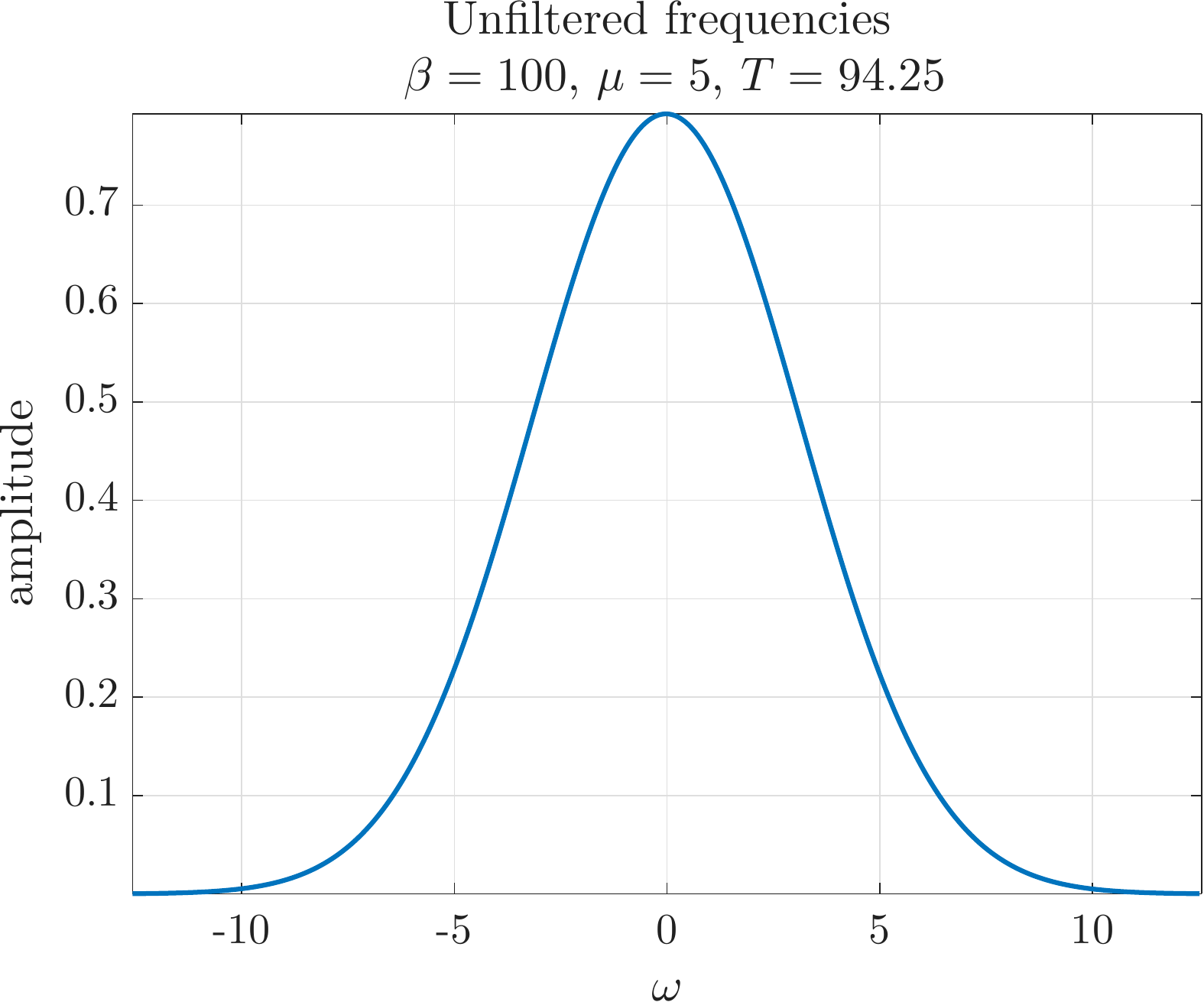}
	\includegraphics[width=7.5cm]{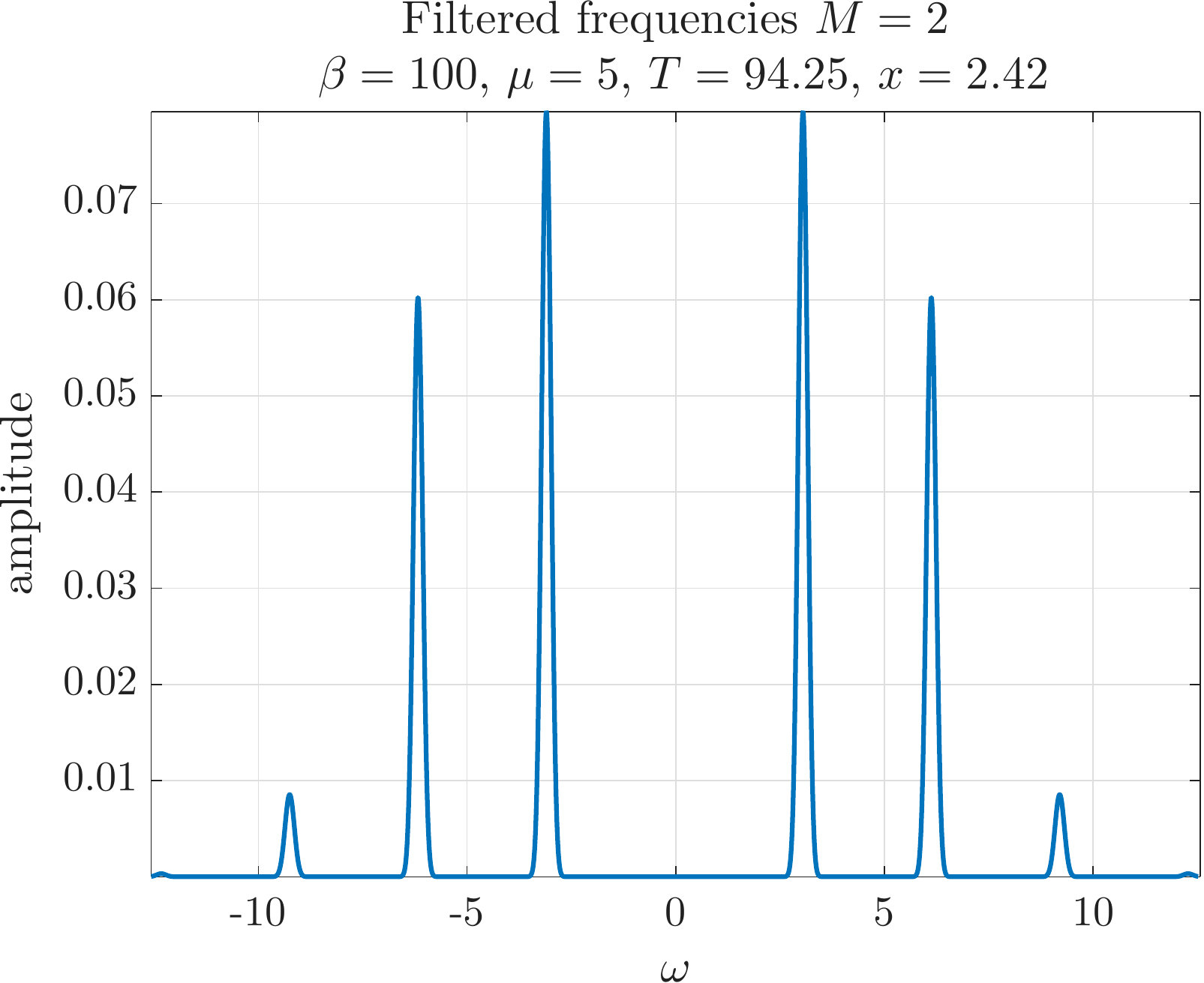}
	\caption{Spectra of the incident wave (left) and the windowed total solution (right) after passing through the $M = 2$ Fabry--Perot filter as in
        Figure \ref{fig:FP_heatmap}.}
	\label{fig:FP_filteration}
\end{figure}

\begin{figure}[th]
	\centering
	\includegraphics[width=7.5cm]{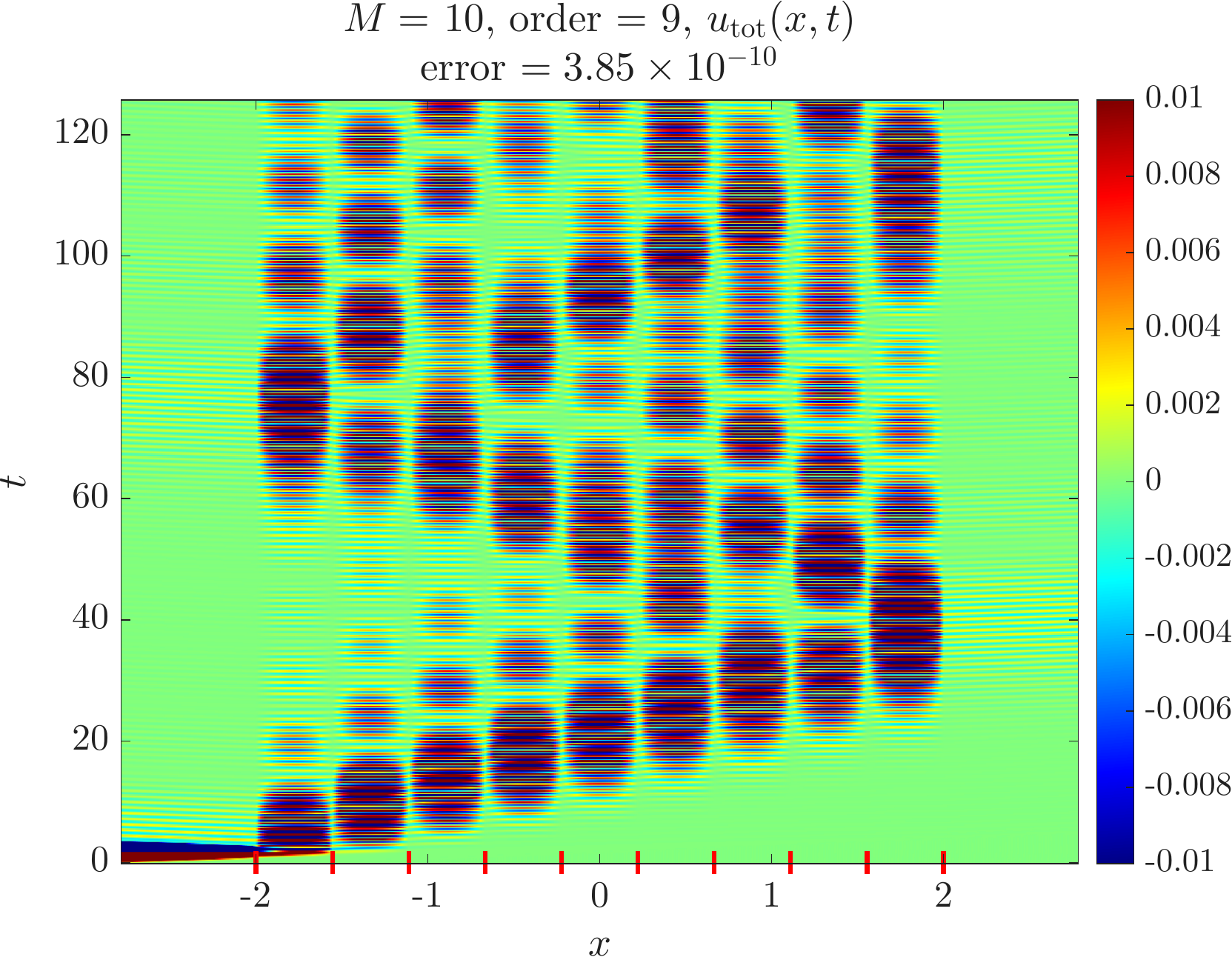}
	\includegraphics[width=7.5cm]{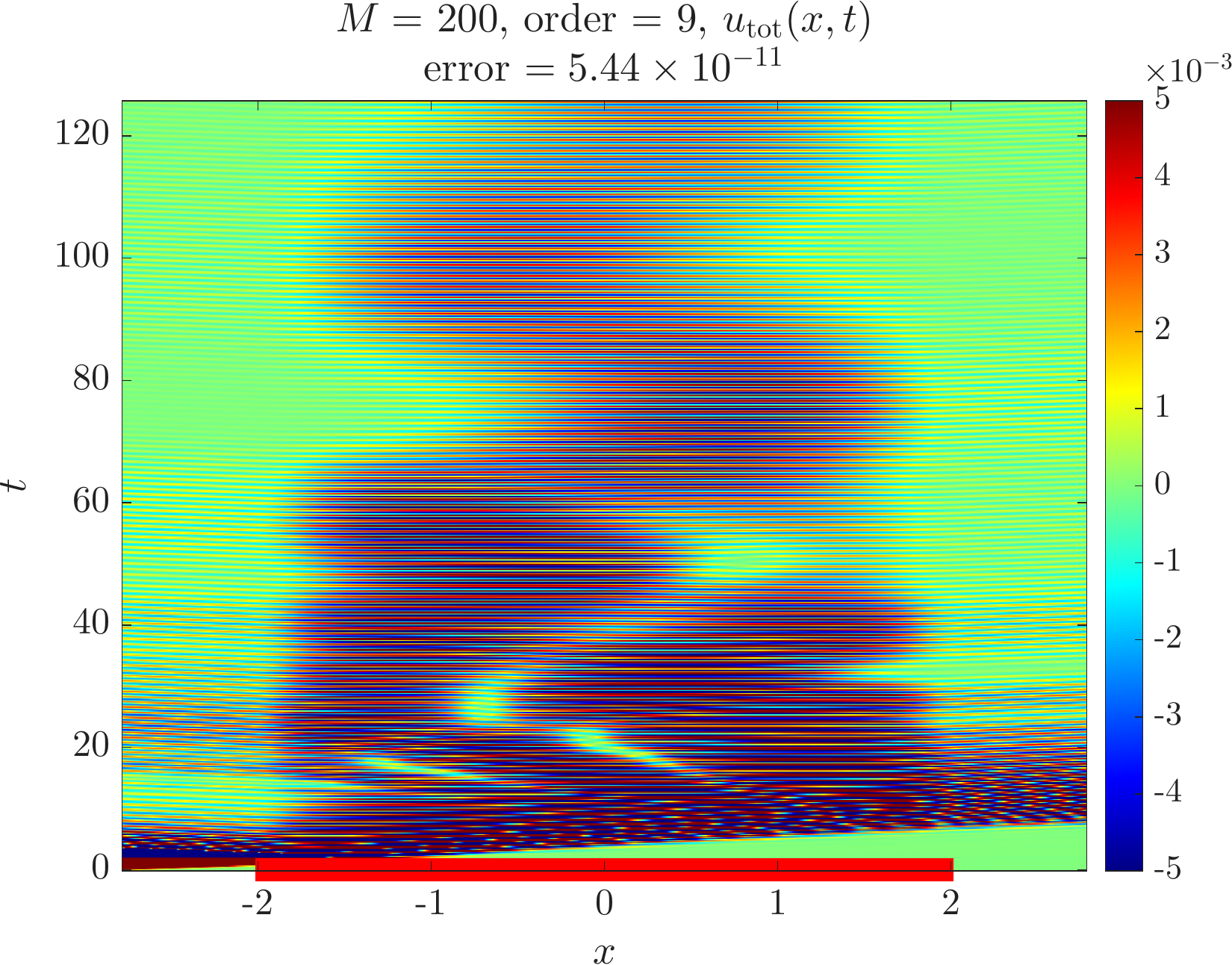}
	\caption{Scattering from $M$ uniformly spaced springs (red tick marks on the $x$-axis) with equal strengths $\beta$.
        The total field is shown, the scheme is 9th order ($p=8$),
        and the estimated accuracies are as shown above the plots.
Left: $M=10$ with $\beta=100$; runtime is 91 seconds.
Right: $M=200$ with $\beta=1$; runtime is 114 seconds.
In both cases 25134 time-steps were used.
}
	\label{fig:FPMS_heatmap}
\end{figure}

\subsection{Scattering and filtering by periodic arrays of springs}\label{sec:uniformSprings}

Regular 1D arrays of identical scatterers are commonly used for
filtering of waves (e.g., the distributed Bragg filter).
Here we apply the WFP scheme to solve scattering problems of this form.
We observe
the transmission filtering, meaning the temporal
Fourier transform of the total wave
$u(x,\cdot)$ at a target point $x$ downstream of the device, as compared
to the Fourier transform of the incident pulse $f$.

The filtering effect depends on the number of sources and the choice of spring constants. For example, given an incident wave of the form~\eqref{eq:incidentPulse} with $\mu = 5$ and $t_0 = -3$, Figure~\ref{fig:FP_heatmap} presents the order 9-accurate total field with two springs located at $x_1 = -0.5$ and $x_2 = 0.5$, and large strengths $\scst_1 = \scst_2 = 100$.
This arrangement of two highly-reflective (near-Dirichlet) scatterers
is common in optics, and is known as a Fabry--Perot interferometer.
Waves get trapped between the two springs; only waves with frequencies
$\om$ close to multiples of $\pi/L$ are transmitted,
$L=x_2-x_1 = 1$ being the cavity length.
Figure~\ref{fig:FP_filteration} compares the Fourier transforms
(approximated via the windowed FFT) of the incident and transmitted
signals for this device.

Figure~\ref{fig:FPMS_heatmap} illustrates experiments with
larger period spring arrays with an incident pulse of the form~\eqref{eq:incidentPulse} where $\mu = 5$ and $t_0 = -4$.
The estimated accuracies (via self-convergence using $\dt/2$)
are again around 9--10 digits and the run times quite short
(see caption and uniform error figures above the plots).
Interesting cavity-coupling behaviors are seen whose
understanding we leave for future work.

We believe that the examples in this section highlight the ease of accurate
wide-band 1D wave simulations with large numbers of point scatterers.

\section{Conclusions and discussion}
\label{sec:conclusions}

We have presented a fast spectrally-accurate history compression scheme (WFP) for
hyperbolic layer potentials for the wave equation.
It reduces the quadratic costs in the number of sources, and in the number of time steps,
both to linear,
with a new cost quasilinear in the number of wavelengths across the domain.
We applied it to the arbitrary-order accurate
solution of problems where 1D waves impinge on a large number of point-like scatterers
that physically model springs attached to a string, with applications in
disordered acoustics and optical systems.
We showed time-dependent solutions for up to 10000 springs with many accurate digits,
and demonstrated the expected cost scaling per time step for up to a million springs.
The scheme would easily incorporate nonlinear scatterers, unlike
solvers that Fourier transform in time (which are intrinsically linear).
The generalizations to 2D and 3D, where point scatterers are replaced by
surface quadrature nodes, are already in progress---%
although for free-space boundary conditions they involve new techniques---%
and will be reported on soon.

The main analytic results are: 1) a proof that for bandlimited densities
WFP is spectrally convergent with respect to Fourier truncation
(Theorem~\ref{t:trunc}), and 2)
various stability and convergence results for the 1st- and 2nd-order cases for two springs
(Section~\ref{sec:stab}).
One future goal is replacing conditions on the densities in Theorem~\ref{t:trunc}
with conditions on the incident wave in the BVP; due to causality,
this will demand more nuanced spaces than strict bandlimited functions.
It would also be interesting to try to generalize our eigenvalue-based stability proofs
to 2nd-order for well-separated springs, orders $p>2$, and $M>2$ springs.

\appendix

\section{Proof of Theorem~\ref{t:trunc} on Fourier truncation}
\label{a:pf}

We use Fourier transforms as defined in \eqref{FT}.
The Kaiser--Bessel
window derivative $\phi'$ in \eqref{phi} has, applying \eqref{KB}, the Fourier transform
\beq
\hatbump(\om) = \frac{b e^{-i\delta\om/2}}{\sinh b} \sinc \sqrt{ \left( \frac{\delta \om}{2}\right)^2 - b^2}
, \qquad \om\in\R
\label{hatbump}
\eeq
We now fix $k\in\Z$, with $|k|>K$.
The Fourier transform of the $k$th influence kernel is related to that of the window
derivative by
\beq
\hat{\Psi}_k(\om) =
\biggl(1-\frac{\om}{2k}\biggr) \hatbump(\om+k) +
\biggl(1+\frac{\om}{2k}\biggr) \hatbump(\om-k),
\qquad k \in \Z \backslash \{0\}, \; \om \in \R,
\label{hatpsigeneral}
\eeq
a formula easily derived from \eqref{psidef} and using the Fourier transform of a
derivative.
Note that the multiplication by sinusoids in \eqref{psidef} shifted the center
frequency of the window from zero to $\pm k$.
In particular, for the Kaiser--Bessel window inserting \eqref{hatbump} gives
\beq
\hat{\Psi}_k(\om) = \frac{b e^{-i\delta\om/2}}{\sinh b}
\left[
 \left(1-\frac{\om}{2k}\right) \sinc \sqrt{(\delta(\om+k)/2)^2 - b^2} +
\left(1+\frac{\om}{2k}\right) \sinc \sqrt{(\delta(\om-k)/2)^2 - b^2}
  \right]
.
\label{hatpsi}
\eeq
Now, by the hypothesis $|k|>K\ge K_0 + 2b/\delta$, the signal frequencies $\om \in [-K_0,K_0]$
lie outside of both shifted window cutoff bands $[k-2b/\delta,k+2b/\delta]$
and $[-k-2b/\delta,-k+2b/\delta]$;
that is,
the arguments of the sinc functions are real-valued, so we may apply $|\sinc x |\le 1/x$.
Thus, doubling the larger of the two $\pm k$ terms, and noting $1+\om/2k < 3/2$,
\beq
\max_{-K_0\le \om \le K_0} \left|\hat{\Psi}_k(\om)\right|
\le
\frac{2b}{\sinh b} \frac{3/2}{\sqrt{(\delta(|k|-K_0)/2)^2 - b^2}}
.
\label{hatpsibnd}
\eeq
Now consider the driving function $F_k$ given by \eqref{Fk}.
The convolution theorem gives $\hat{F}_k(\om) = \hat{\Psi}_k(\om) \hat{S}_k(\om)$,
and, since ${S}_k$ has bandlimit $K_0$,
\beq
\|\hat{F}_k\|_\LLR \le \max_{-K_0\le \om \le K_0} \bigl|\hat{\Psi}_k(\om)\bigr| \cdot
  \|\hat{S}_k\|_\LLR
  \le
  \frac{b}{\sinh b} \frac{3\sqrt{2\pi}MC}{\sqrt{(\delta(|k|-K_0)/2)^2 - b^2}}
  \label{Fkbnd}
\eeq
where in the 2nd step we applied the following bound using Parseval's theorem,
the definition \eqref{densco}, Cauchy--Schwarz,
and the hypothesis that the densities have bounded norm,
\[
\frac{1}{2\pi}\|\hat{S}_k\|^2_\LLR = \| {S}_k\|^2_\LLR \le
\int_\R \left| \sum_{j=1}^M \sigma_j(t) e^{-ikx_j} \right|^2 dt
\le
M \int_\R \sum_{j=1}^M |\sigma_j(t)|^2 dt
= M \sum_{j=1}^M \|\sigma_j\|^2_\LLR \le M^2C^2.
\]

Observe that the Duhamel principle for the ODE system \eqref{FOS}
with initial conditions $\alpha_k(0) = \alpha_k'(0)=0$
gives a new formula for $\alpha_k$ as a convolution of $F_k$ with the propagator,
\beq
\alpha_k(t) = \int_0^t \frac{\sin k(t-\tau)}{k} F_k(\tau) d\tau
\label{newalpha}
\eeq
(contrast the formulae \eqref{eq:alphak} and \eqref{Fk}).
Applying Young's convolution inequality \cite[\S4.2]{liebloss} over $0<t<T$,
\beq
\|\alpha_k\|_\LLT \le \left\| \frac{\sin k\cdot}{k} \right\|_{L^1(0,T)}
\|F_k\|_\LLT
\le
\frac{b}{\sinh b} \frac{3MCT}{k\sqrt{(\delta(|k|-K_0)/2)^2 - b^2}}
\label{akbnd}
\eeq
where we used $\| \sin k\cdot \|_{L^1(0,T)} \le T$, Parseval's theorem, and \eqref{Fkbnd}.

We finally now sum over $k$:
applying the triangle inequality in $\LLT$ to \eqref{EK}, using \eqref{akbnd} gives
\beq
\|E_K\|_\LLT \le \sum_{|k|> K} \|\alpha_k\|_\LLT
\le
6MCT \frac{b}{\sinh b}
\cdot\frac{2}{\delta}
\sum_{k=K+1}^\infty \frac{1}{k\sqrt{(k-K_0)^2 - (2b/\delta)^2}}
.
\label{EKbnd}
\eeq
Shifting the index to $y = k-K_0-2b/\delta$ turns this sum into
\[
\sum_{y=K-K_0-2b/\delta+1}^\infty   
\frac{1}{(y+K_0+2b/\delta) \sqrt{(y+ 2b/\delta)^2 - (2b/\delta)^2}}
  \le
\sum_{y=K-K_0-2b/\delta+1}^\infty \frac{1}{y\cdot y}
\le
\frac{\pi^2}{6},
\]
noting that the sum is over an arithmetic sequence of unit spacing starting
no smaller than 1, using the hypothesis \eqref{Kbnd} on $K$.
Substituting the sum bound into \eqref{EKbnd} gives \eqref{Ebnd}, and the proof is complete.

A couple of notes on the proof are worth mentioning:
\bi
\item
  It is tempting to replace \eqref{newalpha} by the original
  definition \eqref{eq:alphak}, but there would be an obstacle
  since the windowed propagator $\phi(\cdot) \sin k\cdot$ is not in $\LLR$.
\item
  The algebraic prefactors can undoubtedly be improved, for instance
  by splitting the $k$ sum \eqref{EKbnd} into two parts handled separately.
  This would not, however, improve the asymptotic exponential rate $e^{-b}$.
\ei

\bibliographystyle{elsart-num}
\bibliography{refs2}
 
\end{document}